\documentclass[11pt]{amsart}
\usepackage{amsmath,amsthm, amscd, amssymb, amsfonts,mathtools}
\usepackage[all]{xy}
\usepackage{enumitem}

\usepackage{hyperref}
\usepackage{t1enc}
\usepackage[latin1]{inputenc}

\usepackage{fontsmpl}
\usepackage{srcltx}

\allowdisplaybreaks

\newcommand{\toba}{\mathfrak{B}}

\usepackage{color}

\numberwithin{equation}{section}\theoremstyle{plain}

\newtheorem{theorem}{Theorem}[section]

\newtheorem{lem}[theorem]{Lemma}

\newtheorem{pro}[theorem]{Proposition}

\theoremstyle{definition}

\newtheorem{exa}[theorem]{Example}

\theoremstyle{remark}

\newcommand{\ydh}{{}^{H}_{H}\mathcal{YD}}

\newcommand\ev{\operatorname{ev}}
\newcommand\id{\operatorname{id}}

\newcommand\End{\operatorname{End}}

\newcommand\GK{\operatorname{GK-dim}}

\def\ku{\Bbbk}

\def\ot{\otimes}

\def\s{\mathbb{S}}

\def\Z{\mathbb{Z}}
\def\N{\mathbb{N}}

\def\bB{\mathfrak{B}}
\def\hit{\mathfrak{R}}

\def\bX{\mathbf{X}}
\def\bY{\mathbf{Y}}

\newcommand{\J}{{\mathcal J}}

\newcommand\I{\mathbb I}

\def\lg{\langle}
\def\rg{\rangle}

\def\pf{\begin{proof}}
\def\epf{\end{proof}}

\begin{document}


 \title[Nichols algebras associated to solutions  of the  QYBE in rank 3]{Examples of Nichols algebras associated to upper triangular solutions of the Yang--Baxter equation in rank 3}
\author[Giraldi and Silva]{Jo\~ao Matheus Jury Giraldi and Leonardo Duarte Silva}


\address{IME, 
Universidade Federal do Rio Grande do Sul,
Rio Grande do Sul, Brazil.}

\email{joaomjg@gmail.com, leonardoufpel@gmail.com}

\thanks{The work was partially supported by CNPq (Brazil)}

\begin{abstract}
We determine some Nichols algebras that admit a non-trivial quadratic relation associated to some families of upper triangular solutions of the Yang--Baxter equation of dimension 3.
\end{abstract}

\maketitle
\section{Introduction}

Let $\ku$ be a field.
Nichols algebras are Hopf algebras in braided tensor categories with very special properties; they appeared first  in the work of Nichols \cite{N} and play  an important role in the classification of (pointed) Hopf algebras \cite{A, AS2}. 

Recall that a \emph{braided vector space}
is a pair $ (V, \, c) $ 
where $ V $ is a vector space and  $c \in GL ( V \otimes V)$, called a \emph{braiding}, satisfies the \emph{braid equation}:
\begin{align}\label{eqn:braid-bis}
(c\otimes \id)(\id \otimes c)(c\otimes \id) =  (\id \otimes c)(c\otimes \id)(\id \otimes c).
\end{align}
We observe that the braid equation \eqref{eqn:braid-bis} is equivalent to the quantum Yang-Baxter equation (QYBE) 
\eqref{eqn:qybe}.
To each braided vector space one attaches a graded connected braided Hopf algebra, generated in degree one and with all primitive elements in this degree,
called the Nichols algebra of $(V,c)$ or simply of $V$ and denoted by $\toba(V)$. 

Despite this simple definition, the actual computation of the defining relations and the dimension of $\toba(V)$ is a hard task, however needed for the classification program of Hopf algebras. A large amount of work on this question is available, specially when $(V,c)$ is of diagonal type \cite{H-Weyl gpd,H-classif RS,A-jems,A-presentation} but not exclusively, see references in the survey \cite{A}.

\medbreak
Assume that $V$ has finite dimension and let $(v_i, v^i)$ be dual bases for $ V $. The braiding $c$ (or even $V$) is \emph{rigid} if the  map $c^{\flat}: V^* \ot V \to V \ot V^*$ given by $ f\ot v \longmapsto \sum_{i}(\ev\ot \id \ot \id)(f\ot c(v\ot v_i)\ot v^i)$ is bijective.

If  $c$ is rigid, then $\toba(V)$ can be realized as a braided Hopf algebra in the category 
of Yetter-Drinfeld modules $\ydh$ over some Hopf algebra $H$. Thus, a new Hopf algebra
$\toba(V)\# H$ is obtained via the process of \textit{bosonization}.
For this reason we restrict our attention to rigid braidings.

\smallbreak

Assume that $\ku$ is algebraically closed field of characteristic zero.
In \cite{Hi1}, Hietarinta classified all solutions of the QYBE in rank 2.
Based on this result, in \cite{AGi}, the authors considered the rigid braidings which are not of diagonal type. They computed the Nichols algebras associated to them assuming that admit at least one quadratic relation.

Later Hietarinta also classified the non-trivial (rigid) upper triangular solutions of the QYBE for dimension 3 in \cite{Hi2}. 
The following question arises naturally: to study the Nichols algebras related to such solutions. 
It turns out that this is considerably more difficult than the case of rank 2.

In the present paper, we analyze the Nichols algebras associated to some of  the families in \cite{Hi2},
namely $ (\hit_{1, j})_{j\in \I_{10}} $.
As in  \cite{AGi}, we restrict to those that have at least one quadratic relation.
Technically, there is a quadratic relation if and only if $ t^2=1 $ where $t$ is a parameter to be explained in the main part of the paper.
We compute most of the Nichols algebras like this. Here is our main result:

\begin{theorem}\label{th:main}
Let $(V,c)$ be the rigid braided vector space corresponding to some $ R $-matrix $ \hit_{1, j} $. 
Then, $\toba(V)$ has quadratic relations if and only if $ t^2=1 $.
Furthermore, the explicit presentation of $ \toba(V)$, a PBW-basis, the dimension and the GK-dimension are given if $ t^2=1 $ and some extra conditions hold $ ( $when needed$ ) $, see Table \ref{tab:general}.
\begin{table}[ht]
	\caption{Examples of Nichols algebras of rank 3 }\label{tab:general}
	\begin{center}
		\begin{tabular}{| c | c | p{4.8cm} | p{3.7cm} |} 
			Name	& Ref. &  $\toba(V)$   &   Extra Conditions  	\\\hline
			\hline
			$\hit_{1, 1}$  & \S \ref{subsec:hit11} & Table \ref{tab:h11} &   \\ \hline
			$\hit_{1, 2}$   & \S \ref{subsec:hit12} & Propositions \ref{prop:case1,2_azul}, \ref{prop:case1,2_vermelho}, Table \ref{tab:h12_3} &  
			\\ \hline
			$\hit_{1, 3}$  
			& \S \ref{subsec:hit13} & Table \ref{tab:h13} &  \\		
			\hline
			$\hit_{1, 4}$  & \S \ref{subsec:hit14} & Table \ref{tab:h14} & 
			\\ \hline
			$\hit_{1, 5}$   & \S \ref{subsec:hit15} & Propositions \ref{prop:case1,5_t=1_and_p=2q}, \ref{prop:case1,5_t=-1_and_p_neq_0}  
			& $ t=1 $ and $ p=2q $, or  $ t=-1 $ and $ p=0 $ 
			\\ \hline
			$\hit_{1, 6}$   & \S \ref{subsec:hit16} & Table \ref{tab:h16} & 
			\\ \hline
			$\hit_{1, 7}$ 
			& \S \ref{subsec:hit17} & Table \ref{tab:h17}& $ t=1 $ and $ a=p-2q $,  or $ t=-1 $ and $ a=-p $ \\\hline 
			$\hit_{1, 8}$   & \S \ref{subsec:hit18} & Table \ref{tab:h18} & 
			\\ \hline 
			$\hit_{1, 9}$   & \S \ref{subsec:hit19} & Table \ref{tab:h19} & $ b=1 $ 
			\\ \hline 
			$\hit_{1, 10}$   & \S \ref{subsec:hit110} & Table \ref{tab:h110} &  
			\\ \hline
			\end{tabular}
			\end{center}
			\end{table}
\end{theorem}

To determine these Nichols algebras, we use the following technique, already used in many papers.

First, we look for relations in small degrees. For this, we use the quantum symmetrizers $ Q_n $ \eqref{eqn:simetrizador} and the following convenient fact:
Let $x \in T(V)$. If $ \partial_{f} (x) =0  $ for all $  f\in V^* $, then $ x \in \J(V) $.
See \S \ref{subsec:Nichols} for details and notation.
Thus, we have a set of homogeneous relations. Write $I$ for the ideal generated by such relations. 

Now, let $\bB = T(V)/ I$ be the pre-Nichols algebra and $\pi: \bB \twoheadrightarrow \toba(V)$ be the natural projection.
Suppose that there is a relation $0 \neq r\in \ker \pi$ which we assume homogeneous of minimal degree. Then we have that $\partial_{f} (r) =0$ by minimality of the degree. If this leads to a contradiction, then $\pi$ is actually an isomorphism.
Otherwise, we get a new relation and repeat the process.
\smallbreak

Finally, we also mention that the constructions made in Subsections \ref{subsubsec:hit12_azul} and \ref{subsubsec:hit12_vermelho} can be generalized for the cases which are excluded by the extra conditions. 
We obtain analogous relations to the ones presented in Propositions \ref{prop:case1,2_azul} and \ref{prop:case1,2_vermelho}.
However, we were not able to prove that they are enough.

\smallbreak

The paper is organized as follows. In Section \ref{sec:preliminaries} we recall definitions and facts that are used throughout the paper.
Then, in Section \ref{sec:quadratic-relations} we do the case-by-case analysis of the Nichols algebras with quadratic relations
of the braided vector spaces $(V,c)$ corresponding to the families $ (\hit_{1, j})_{j\in \I_{10}} $.

\subsection*{Notation}
Let $ j, k\in\Z $.
If $j \le k $, then we denote $\I_{j, k} =\{j, j+1, \dots, k\}$ and $\I_k = \I_{1,k}$.
If $j > k $, then  $\I_{j, k}=\emptyset$.
Sums and products  over an empty set of indices are 0 and 1, respectively.

The maps $ \chi_{_{o}}, \,  \chi_{_{e}}: \mathbb{Z}\longrightarrow \{0, \, 1\}$ stand for the odd and even characteristic functions.
We denote the floor function by $ \left\lfloor \ \right\rfloor: \mathbb{R}\to \mathbb{Z} $.

\subsection*{Acknowledgments} 
The authors thank Nicol\'as Andruskiewitsch and Iv\'an Angiono for fruitful conversations at different moments of our study.

\section{Preliminaries}\label{sec:preliminaries}

\subsection{Yetter-Drinfeld modules and Nichols algebras}\label{subsec:Nichols}

Fix $H$ a Hopf algebra. A \textit{left Yetter-Drinfeld module} $M$ over $H$
is a left $H$-module $(M,\cdot)$ and a left $H$-comodule
$(M,\lambda)$ that satisfies
the compatibility condition 
\begin{align*}
\lambda(h\cdot m) &= h_{(1)}m_{(-1)}\mathcal{S}(h_{(3)})\ot h_{(2)}\cdot m_{(0)},
& & h\in H,\, m\in M. 
\end{align*}
The related category is denoted by $\ydh$.
If $ \mathcal{S} $ is bijective, then $\ydh$ is  braided and monoidal: the braiding $c_{M,N}:M\ot N \to N\ot M$ is given by
$c_{M,N}(m\ot n) = m_{(-1)}\cdot n \ot m_{(0)}$ for all
$m\in M,n\in N$ and $M, N \in \ydh$. 
A Hopf algebra in $\ydh$ is called a \textit{braided} Hopf algebra.

\smallbreak

Given $V \in \ydh$, we say that a braided $\N$-graded Hopf algebra $R = \bigoplus_{n\geq 0} R^{(n)} \in \ydh$  is a \textit{Nichols algebra} for $V$ if  $\Bbbk\simeq R^{(0)}$, $V\simeq R^{(1)} \in \ydh$,
$R^{(1)} = \mathcal{P}(R)$ and 
$R$ is generated as algebra by $R^{(1)}$.

Nichols algebra always exists and is unique up to isomorphism.
It is usually denoted by 
$\toba(V) = \bigoplus_{n\geq 0} \toba^{n}(V) $ and
is given by the quotient of the tensor algebra $T (V)$ by the largest 
homogeneous Hopf ideal  $\J (V)= \bigoplus_{n\geq 2} \mathcal{J}_{n}(V)$  generated by 
homogeneous elements of degree $\geq 2$.

The ideal $\J(V) $ has an alternative description: $\J_n(V)$ is the kernel of the $ n$th quantum symmetrizer associated to $c$ 
\begin{align}\label{eqn:simetrizador}
	Q_{n} &= \sum_{\sigma \in \s_{n}}\rho_{n}(M (\sigma)) \in \End(T^{n}(V)), \qquad\qquad n\geq 2,
\end{align}
where $\rho_{n}:\mathbb{B}_{n} \to GL(T^{n}(V)), \, n \geq 2 $, is the representation induced by $ c $ on the braid group $ \mathbb{B}_{n} $, and $ M $ is the Matsumoto section
corresponding to the canonical projection $\mathbb{B}_{n}\twoheadrightarrow \s_{n}$.
In particular,
\begin{align}\label{eq:J2}
\J_2(V) &= \ker\, (\id + c).
\end{align}

\smallbreak

Another way to find relations in $ \toba(V) $ is through left skew derivations.
Let $f \in V^*$ and consider $ \partial_{f} \in \End  T(V)$ given by
\begin{align*}
\begin{aligned}
\partial_f(1) &=0,& &\partial_f(v) = f(v), \quad v\in V,
\\
\partial_f(xy) &= \partial_f(x)y + \sum x_i \partial_{f_i}(y), & &\text{where } c^{-1} (f \ot x) = \sum x_i \ot f_i.
\end{aligned}\end{align*}

A \textit{pre-Nichols algebra} ${\toba}$ is any graded braided Hopf algebra intermediate between $T(V)$ and $\toba(V)$; in other words, any
braided Hopf algebra of the form $T(V)/I$ where $I \subseteq \mathcal{J}(V)$ is a homogeneous Hopf ideal.

The skew derivations are well behaved with respect to pre-Nichols algebras, that is, for all $ f\in V^* $ we can define $\partial_{f} \in \End  \bB$ satisfying the properties above and compatible with the projection
$T(V) \twoheadrightarrow \bB$. 

Finally, we recall the following result which is very used to determine Nichols algebras. 
Let $x\in \bB^{n} $ with $n\geq 2$. If $\partial_{f}(x)=0$ for all $f\in V^{*}$, then  $x\in \mathcal{J}_{n}(V)$.

See  \cite{A}, \cite{AS1} for more details.

\subsection{Upper triangular solutions of the Yang--Baxter equation in rank 3} \label{subsec:solutions}

Let $R \in \End(V \otimes V)$. 
We say that the map $ R $ satisfies the quantum Yang-Baxter equation (QYBE) if
\begin{align}\label{eqn:qybe}
R_{12}R_{13}R_{23} &=  R_{23}R_{13}R_{12}.
\end{align}
A solution of the QYBE is also called a $ R $-matrix.
The QYBE \eqref{eqn:qybe} is equivalent to the braid equation \eqref{eqn:braid-bis} via
\begin{align*}
R \longleftrightarrow c &= \tau R,
\end{align*}
where $\tau: V\ot V \to V \ot V$ is the usual flip, $\tau(x \ot y) = y \ot x$.


In \cite{Hi2}, Hietarinta classified all the non-trivial upper triangular solutions of the QYBE in rank 3.
It turns out that there are 35 families of solutions, enumerated by $(\hit_{1, j})_{j\in \I_{10}}$, $(\hit_{2, j})_{j\in \I_{9}}$, $(\hit_{3, j})_{j\in \I_4}$, $(\hit_{4, j})_{j\in \I_4}$, $(\hit_{6, j})_{j\in \I_2}$, $(\hit_{7, j})_{j\in \I_2}$, $(\hit_{8, j})_{j\in \I_3}$ and $\hit_{9, 1}$.
All of them are invertible and rigid.

For this study, we focus on the  families $(\hit_{1, j})_{j\in \I_{10}}$. 
Notice that we homogenize  the original solutions by a parameter $ t $.

Given a solution of the Yang-Baxter equation 
$$ \hit= \left(\begin{array}{ccccccccc}
R^{11} & R^{21} & R^{31} & R^{12} & R^{22} & R^{32} & R^{13} & R^{23} & R^{33} 
\end{array} \right),$$ 
where each column matrix $ R^{ij} $ is the transpose of 
$$ \left(\begin{array}{ccccccccc}
r^{ij}_{11} & r^{ij}_{12} & r^{ij}_{13} & r^{ij}_{21} & r^{ij}_{22} & r^{ij}_{23} & r^{ij}_{31} & r^{ij}_{32} & r^{ij}_{33} 
\end{array} \right), $$ 
we associate the following braiding 
\begin{align} \label{eqn:QYBE_to_braiding}
c(x_i \ot x_j) = \sum_{\mathclap{k, \, \ell \,\in\, \I_3}} r^{ij}_{k\ell} \,  x_k\ot x_\ell . 
\end{align}

\subsection{Colexicographic order}\label{subsec:colex_order}
Let $ (A,\prec_A) $ and $ (B,\prec_B) $ be two partially ordered sets.
Define on the Cartesian product $ A\times B $ the following relation:
\begin{align*}
(a_1,b_1)\prec_{A\times B}^{colex}(a_2,b_2) \text{ iff }
b_1\prec_{B}b_2\text{ or } (b_1=b_2 \text{ and } a_1\prec_{A}a_2).
\end{align*}
It is a partial order and known as the \textit{colexicographical order} on $ A\times B $.
If $ A $ and $ B $ are totally ordered, then it is a total order also.

More generally, one can define the colexicographic order on set of all (non-commutative) finite words $A^{\N}$ with alphabet $ A $.
Write $ \ku A^{\N} $ for the corresponding algebra of finite words.

Let $ x\in \ku A^{\N}, \, x\neq 0 $.
Then there is a unique writing $ x = \sum_{i\in \I_n} \lambda_i \, x_i,\, \lambda_i \in\ku,\,  x_i\in A^\N $, such that $ \lambda_{n}\neq 0 $ and $ x_j\prec x_n $ for all $ j\in \I_{n-1} $. 
We call $ x_n $ the maximal term of $ x $ and denote it by $ \max x $.

Observe that the same constructions can be done for the context of commutative words/polynomials.

\begin{exa} Let $ A=\{a\prec b\} $ be a totally ordered set. 
With respect to the colexicographic order, we have that
\begin{align*}
1 \prec\! a\prec \! a^2 \prec\! a^3 \prec\! ba^2 \prec\! ba\prec\! aba \prec\! b^2a \prec\! b\prec\! ab\prec\! a^2b \prec\! bab \prec\! b^2\prec\! ab^2 \prec\! b^3\! ,
\end{align*}
where $ 1 $ stands for the word of length $ 0 $.
\end{exa}

\section{Nichols algebras of rank three with quadratic relations}\label{sec:quadratic-relations}

For this section, we set the following notation:
\begin{align}\label{eqn:B_0}
B_0 &=\{x_1^{a_1}x_2^{a_2}x_3^{a_3}: \, 0 \leq a_i < 2\}; \\ \label{eqn:B_1} 
B_1 &=\{x_1^{a_1}x_2^{a_2}x_3^{a_3}: \, a_i\geq 0; \,  a_1, a_2 < 2 \};\\ \label{eqn:B_3}
B_3 &=\{x_1^{a_1}x_2^{a_2}x_3^{a_3}: \,  a_i \geq 0 \}.
\end{align}

We always assume that the homogenizing parameter $ t\neq 0 $.

\subsection{Case $ \hit_{1, 1}$}\label{subsec:hit11}
Let $ c $ be the braiding associated to the solution of QYBE
\begin{align*}
	\hit_{1, 1} = \left(\begin{array}{ccc|ccc|ccc}
	t & \cdot & ta & \cdot & \cdot & \cdot & -ta & \cdot & -tab \\ 
	\cdot & t & \cdot & \cdot & \cdot & t(a-b) & \cdot & -ta & tp \\ 
	\cdot & \cdot & t & \cdot & \cdot & \cdot & \cdot & \cdot & -tb  \\ \hline
	\cdot & \cdot & \cdot & t & \cdot & ta & \cdot & t(b-a) & -tp \\ 
	\cdot & \cdot & \cdot & \cdot & t & \cdot & \cdot & \cdot & \cdot \\ 
	\cdot & \cdot & \cdot & \cdot & \cdot & t & \cdot & \cdot & \cdot \\ \hline
	\cdot & \cdot & \cdot & \cdot & \cdot & \cdot & t & \cdot & tb  \\ 
	\cdot & \cdot & \cdot & \cdot & \cdot & \cdot & \cdot & t & \cdot  \\ 
	\cdot & \cdot & \cdot & \cdot & \cdot & \cdot & \cdot & \cdot & t 
	\end{array} \right).
\end{align*}  
 Just for this initial case, we write the braiding explicitly. See \eqref{eqn:QYBE_to_braiding} for details.
\begin{align*}
c(x_1 \ot x_1) &= t \, x_1 \ot x_1; \qquad\qquad\qquad\qquad\qquad\quad \ \,  c(x_2 \ot x_1) = t \, x_1 \ot x_2; \\
c(x_3 \ot x_1) &= t \, x_1 \ot x_3 + ta \, x_1\ot x_1; \qquad\qquad\qquad c(x_1 \ot x_2) = t \, x_2 \ot x_1; \\
c(x_2 \ot x_2) &= t \, x_2 \ot x_2; \\
c(x_3 \ot x_2) &= t \, x_2 \ot x_3 + ta \, x_2\ot x_1 + t(a-b) \, x_1\ot x_2; \\ 
c(x_1 \ot x_3) &= t \, x_3 \ot x_1 - ta \, x_1\ot x_1; \\
c(x_2 \ot x_3) &= t \, x_3 \ot x_2 + t(b-a) \, x_2\ot x_1 -ta \, x_1\ot x_2; \\
c(x_3 \ot x_3) &= t \, x_3 \ot x_3 + tb \, x_3\ot x_1 -tp \, x_2\ot x_1 - tb \, x_1\ot x_3 \\ 
& \ \ \ +tp \, x_1\ot x_2 - tab \, x_1\ot x_1.
\end{align*}

\begin{pro} \label{prop:case1,1}
If $ t^2 \neq 1 $, then there are no quadratic relations. Otherwise, the Nichols algebras are as in Table \ref{tab:h11}, where 
\begin{align}\label{eqn:Nichols1,1_t=-1}
\lg  x_2x_1 + x_1x_2, \, x_3x_1 + x_1x_3, \, x_3x_2 \, +  & \,x_2x_3 - b \, x_1x_2,  \\ 
\notag &x_1^2,\, x_2^2, \, x_3^2 -b\, x_1x_3 + p\, x_1x_2  \rg;  \\ \label{eqn:Nichols1,1_t=1}
\lg  x_2x_1 - x_1x_2, \, x_3x_1 - x_1x_3 - a\, x_1^2,\,  &\, x_3x_2 - x_2x_3 + (b-2a) \, x_1x_2\rg.
\end{align}	

\begin{table}[ht]
	\caption{Nichols algebras of type $ \hit_{1, 1.} $}\label{tab:h11}
	\begin{center}
		\begin{tabular}{| c | c | c | c | c | c |}\hline
			$ t $  & $\J(V)$ &  Basis &$\GK$
			\\\hline
			$-1$ &  \eqref{eqn:Nichols1,1_t=-1} & $ B_0 $ \eqref{eqn:B_0} & $ 0, \, \dim = 8 $ \\\hline
			$1$ &  \eqref{eqn:Nichols1,1_t=1} & $ B_3  $ \eqref{eqn:B_3} & $ 3 $ \\\hline
		\end{tabular}
	\end{center}
\end{table}		
\end{pro}


\pf The quadratic relations are obtained by \eqref{eq:J2}. 
If $ t=-1 $, then the pre-Nichols algebra $ T(V)/I $ where I is the ideal \eqref{eqn:Nichols1,1_t=-1} is a 8-dimensional algebra whose basis is $ B_0 $. In this case, it is clear that no more relations exist and then $ \toba (V) = T(V)/I $.

If $ t=1 $, then we write $ \bB = T(V)/I $ where I is the ideal \eqref{eqn:Nichols1,1_t=1}.
Observe that $ B_3 $ generates linearly $ \bB $ since the relations
\begin{align}\label{eqn:relations_prop_1.1}
	x_3x_1^{a_1} = x_1^{a_1}x_3 + a_1a \, x_1^{a_1+1}  \text { and } \  
	x_3x_2^{a_2} = x_2^{a_2}x_3 + a_2(2a-b) \, x_1x_2^{a_2} 
\end{align}
hold in $ \bB $, for all $ a_1, a_2\geq 0 $.
An easy calculation shows that
\begin{align*}
c(x_1^{a_1}x_2^{a_2}\ot x_3) &= (x_3- a(a_1+a_2)x_1)\ot x_1^{a_1}x_2^{a_2} + a_2(b-a) \, x_2\ot x_1^{a_1+1}x_2^{a_2-1},
\end{align*}
for $ a_1, a_2\geq 0 $, what implies that
\begin{align}\label{eqn:partial_prop_1.1}
\partial_3(x_1^{a_1}x_2^{a_2}x_3^{a_3}) &= x_1^{a_1}x_2^{a_2}\partial_3(x_3^{a_3}), \qquad a_1, a_2, a_3\geq 0.
\end{align}
As $ \partial_3(x_3^{a_3}) \in \bB_{a_3-1}, $ there are $ \mu_{a_3, c_1, c_2, c_3}\in \ku $, with $ c_1+c_2+c_3=a_3-1 $ and $ c_i\geq 0 $, such that
\begin{align}\label{eqn:sum_prop_1.1}
\partial_3(x_3^{a_3}) &= \sum_{\mathclap{\substack{c_1+c_2+c_3=a_3-1, \, c_i\geq 0 }}}  \mu_{a_3, c_1, c_2, c_3}\, x_1^{c_1}x_2^{c_2}x_3^{c_3}, \qquad a_3 \geq 1.
\end{align}
We claim that $ \mu_{a_3, 0, 0, a_3-1} = a_3 $.
Indeed, the case $ a_3=1 $ is obvious. Assume that $ \mu_{a_3, 0, 0, a_3-1} = a_3 $, then
\begin{align*}
\partial_3(x_3^{a_3+1}) &= x_3^{a_3} + (x_3+ b\, x_1)\partial_3(x_3^{a_3}) \\
&= x_3^{a_3} + (x_3+ b\, x_1) \sum_{\mathclap{\substack{c_1+c_2+c_3=a_3-1, \, c_i\geq 0}}} \mu_{a_3, c_1, c_2, c_3}\, x_1^{c_1}x_2^{c_2}x_3^{c_3}.
\end{align*}
By \eqref{eqn:relations_prop_1.1}, the claim holds.

Suppose that the natural projection  $ \pi: \bB \to \toba (V) $ is not an isomorphism. Pick a linear homogeneous relation  of minimal degree $ n\geq 3 $
\begin{align*}
r =  \sum_{\mathclap{\substack{a_1+a_2+a_3=n, \, a_i\geq 0}}}  \lambda_{a_1, a_2, a_3}\, x_1^{a_1}x_2^{a_2}x_3^{a_3}.
\end{align*}
Thus, by \eqref{eqn:partial_prop_1.1} and \eqref{eqn:sum_prop_1.1},
\begin{align}\label{eqn:partial_sum_prop_1.1}
 0 &= \partial_3(r) = 
 \sum_{ \mathclap{\substack{a_1+a_2+a_3=n \\ a_i\geq 0}}}   \lambda_{a_1, a_2, a_3}   \ \, \sum_{\mathclap{\substack{c_1+c_2+c_3=a_3-1 \\ c_i\geq 0}}} \mu_{a_3, c_1, c_2, c_3} \, x_1^{a_1+c_1}x_2^{a_2+c_2}  x_3^{c_3}.
\end{align}

Observe that the term $ x_3^{n-1} $ appears only one time in \eqref{eqn:partial_sum_prop_1.1}. 
By the minimality of $ n $ and the claim above, we get $ \lambda_{0, 0, n}=0 $. 
Inductively on $ k $, we prove that $ \lambda_{j, k-j, n-k}=0,\, j\in\I_{0, k} $, for all $ k= 0, 1, \cdots, n-1 $ since, at each step $ k $, the term $ x_1^{j}x_2^{k-j} x_3^{n-k} $ shows up just one time in \eqref{eqn:partial_sum_prop_1.1}.

Hence, the relation $ r $ can be rewritten as $ r =  \sum
\lambda_{a_1, a_2, 0}\, x_1^{a_1}x_2^{a_2} $.
Since $ \partial_1(x_1^{a_1}x_2^{a_2}) = a_1 \,x_1^{a_1 -1}x_2^{a_2} $ and $ \partial_2(x_1^{a_1}x_2^{a_2}) = a_2 \, x_1^{a_1}x_2^{a_2-1} $, we conclude that $ r=0 $, a contradiction.
Therefore, there are no more relations and $ \pi $ is an isomorphism.
\epf

\subsection{Case $ \hit_{1, 2}$}\label{subsec:hit12}
Let $ c $ be the braiding associated to the solution of QYBE  
$$ \hit_{1, 2} = \left(\begin{array}{ccc|ccc|ccc}
t & \cdot & tb & \cdot & \cdot & \cdot & tp & \cdot & ta \\ 
\cdot & t & \cdot & \cdot & \cdot & t(p-q) & \cdot & tq & tk \\ 
\cdot & \cdot & t & \cdot & \cdot & \cdot & \cdot & \cdot & tp  \\ \hline
\cdot & \cdot & \cdot & t & \cdot & tb & \cdot & \cdot & -tk \\ 
\cdot & \cdot & \cdot & \cdot & t & \cdot & \cdot & \cdot & \cdot \\ 
\cdot & \cdot & \cdot & \cdot & \cdot & t & \cdot & \cdot & \cdot \\ \hline
\cdot & \cdot & \cdot & \cdot & \cdot & \cdot & t & \cdot & tb  \\ 
\cdot & \cdot & \cdot & \cdot & \cdot & \cdot & \cdot & t & \cdot  \\ 
\cdot & \cdot & \cdot & \cdot & \cdot & \cdot & \cdot & \cdot & t 
\end{array} \right). $$
For the explicit presentation of $ c $, use \eqref{eqn:QYBE_to_braiding}.

\smallbreak


By \eqref{eq:J2}, it is clear that the Nichols algebra associated to $ c $ has a quadratic relation if and only if $ t^2=1 $.
To simplify the study of these Nichols algebras, we consider three subcases:
\begin{enumerate}[leftmargin=*,label=\rm{(\roman*)}] 
	\item\label{subcase:1,2_1} $ t= 1 $ and $ b\neq  p-2q;$
\item\label{subcase:1,2_2} $ t= -1 $ and $ b\neq  -p;$
\item\label{subcase:1,2_3} $ t= 1 $ and $ b=  p-2q$, or $ t= -1 $ and $ b=  -p$.
\end{enumerate}
Each subcase is treated separately in the next three subsections. 

\subsubsection{Case \ref{subcase:1,2_1}}\label{subsubsec:hit12_azul}
Assume that
$ t=1 $ and $ b\neq  p-2q$.
Here, we classify the Nichols algebras under these conditions.
To do so, first we consider a suitable environment  to deal with the problem. 
Then, we show two lemmas needed to prove the main result which is Proposition \ref{prop:case1,2_azul}.

\smallbreak
By \eqref{eq:J2}, we get that 
\begin{align}\label{eqn:1.2_1_relation degree 2}
&x_2x_1 - x_1x_2,  &&x_3x_1 - x_1x_3 +\frac{(p-b)}{2}x_1^2 
\end{align}
are the quadratic relations of $ \bB(V) $. For $ n\geq 1 $, define recursively
\begin{align*}
z_n= (x_3-nb\, x_1)z_{n-1} -z_{n-1}x_3 +(n-1)!(q-p)\left( \frac{-b-p}{2}\right)^{n-1} \! \! x_1^nx_2,
\end{align*}
where $ z_0:=x_2 $. 
By induction on $ n $, we prove that the terms $ z_n $ satisfy the following properties in $ \bB(V) $:
\begin{align}\label{eqn:1.2_1_properties1}
\begin{split}
z_nx_i = x_iz_n, \qquad\qquad
z_nz_m = z_nz_m, \ \ \ \ \ \ \ \ \ \ \ \ \ \ \ \ \ \ \ \ \ \\ 
c(z_n\otimes x_i) = x_i\otimes z_n, \quad\quad
c(z_n\otimes x_3) = (x_3+(np+q)x_1)\otimes z_n, \\ 
\partial_3(z_n) = 0, \qquad\quad
\partial_2(z_n) = \gamma_n\, x_1^n,\qquad\quad
\partial_1(z_n) = \beta_n\, z_{n-1},  \,\,\,\,\,\,
\end{split}
\end{align}
for all $ n\geq 1 $, $ m\in\I_{n-1} $ and $ i\in \I_2 $, where
\begin{align}\label{eqn:1.2_1_gamma_and_pi}
\gamma_n= n!(q-p)\left( \frac{-p-b}{2}\right) ^{n-1} \!\!\!\!\!\!\text{ and } \, \beta_n= -n\left( \frac{n+1}{2}\, b+\frac{n-1}{2}\, p+q\right).
\end{align}

We leave to reader the proof of these properties.
We just observe that 
\begin{align}\label{eqn:1.2_1_relation zn zn-1}
	z_{n+1}z_n &- z_nz_{n+1}, &&n\geq 0,
\end{align}
are new relations and proved through derivations.
We also mention some additional identities that are necessary to prove the properties \eqref{eqn:1.2_1_properties1}:
\begin{align} \label{eqn:1.2_1_ideal a esquerda1}
\begin{split}
	c(x_1^n\otimes x_3) &= (x_3+np\, x_1)\otimes x_1^n, \\
	x_3x_1^n &= x_1^nx_3 -n\left(\frac{p-b}{2}\right) x_1^{n+1},
\end{split}
&&n\geq 0.
\end{align}

Set $ \bB = T(V)/I$ where $ I\subseteq\J(V) $ is the ideal generated by \eqref{eqn:1.2_1_relation degree 2} and \eqref{eqn:1.2_1_relation zn zn-1}.
Observe that properties \eqref{eqn:1.2_1_properties1} remain true in $ \bB $.
In particular,
\begin{align}\label{eqn:1.2_1_ideal a esquerda2}
\begin{split}
	x_3z_n^m =  \, z_n^mx_3 + m\, z_n^{m-1}z_{n+1}+ m(n+1)b \, x_1z_n^m \\
	-m\, n!(q-p)\left( \frac{-b-p}{2}\right)^n x_1^{n+1}x_2z_n^{m-1}
\end{split}	
\end{align}
holds in  $ \bB $ for all $ n, m \geq 0 $.
Hence, it follows that 
\begin{align*}
 B_\infty =\{x_1^{a_1}x_2^{a_2}z_1^{b_1}z_2^{b_2} \cdots z_n^{b_n}x_3^{a_3}: n\geq 1; \, a_i, b_j\geq 0 \}
\end{align*}
generates linearly the pre-Nichols algebra $ \bB $ since by \eqref{eqn:1.2_1_ideal a esquerda1} and \eqref{eqn:1.2_1_ideal a esquerda2}, $ \Bbbk B_\infty $ is a left ideal of $ \bB $ which contains its unit.

Using the properties \eqref{eqn:1.2_1_properties1}, we prove that
\begin{align*}
\begin{split}
\partial_i(z_n^{b_n})&=  b_n \,\partial_i(z_n) z_n^{b_n-1}, \qquad  i\in\I_3, \, n\geq 1. 
\end{split} 
\end{align*}
Thus, we obtain that
\begin{align*}
\begin{split}
\partial_i(x_1^{a_1}x_2^{a_2}z_1^{b_1} \cdots  z_n^{b_n}) = \  & \partial_i(x_1^{a_1})x_2^{a_2}z_1^{b_1} \cdots z_n^{b_n}+ x_1^{a_1}\partial_i(x_2^{a_2})z_1^{b_1} \cdots z_n^{b_n} \\
&+ x_1^{a_1}x_2^{a_2}\partial_i(z_1^{b_1}) \cdots z_n^{b_n}+ \cdots + x_1^{a_1}x_2^{a_2}z_1^{b_1} \cdots \partial_i(z_n^{b_n}),
\end{split}
\end{align*} 
for all $ i\in\I_3 $. In particular, $ \partial_3(x_1^{a_1}x_2^{a_2}z_1^{b_1} \cdots z_n^{b_n}) = 0 $.
 
\begin{lem}\label{lem:1.2_1_partial 3} For $ a_i, b_j \geq 0 $ and $ n>0 $, the following identities hold in $ \bB $: 
	\begin{align}\label{eqn:1_lem:1.2_1_partial 3}
	\partial_i(\bY x_3^{a_3}) &= \partial_i(\bY)\, x_3^{a_3} + \bY(\partial_i(x_3^{a_3}) +\delta_{i, 1}\alpha\, \partial_3(x_3^{a_3}) ),  \qquad\qquad\, i\in\I_3,
 \\ \label{eqn:4_lem:1.2_1_partial 3}
	\partial_3(x_3^{n}) &=    \sum_{\mathclap{i\in\I_{0, n-1}}}\binom{n}{i+1} x_1^i x_3^{n-1-i} \prod_{j\in \I_{i}} \left(  b+(j-1)\, \frac{b-p}{2}\right),
	\end{align} 
where $ \bY = x_1^{a_1}x_2^{a_2}z_1^{b_1} \cdots z_n^{b_n} $ and $ \alpha = a_1p+a_2q+\sum_{j\in \I_n} b_j(jp+q) $. 
\end{lem} 
 
\pf 
For all $ n, m\geq 0 $, we have
\begin{align*}
	c(z_n^m\otimes x_3) &= (x_3+m(np+q)x_1)\otimes z_n^m,  \\
	c(x_1^m\otimes x_3) &= (x_3+mp\, x_1)\otimes x_1^m.
\end{align*}
Then, $
c(\bY\otimes x_3) = (x_3+ \alpha x_1)\otimes \bY .
$
In particular, \eqref{eqn:1_lem:1.2_1_partial 3} follows 
since
$
c(\bY\otimes x_i) = x_i\otimes \bY,  $ for each $i\in\I_2.
$

To obtain \eqref{eqn:4_lem:1.2_1_partial 3}, first observe that $ \partial_3(x_3^{n+1}) = (bx_1+x_3)\partial_3(x_3^{n}) $. Then, 
%
\begin{align*}
\partial_3(x_3^{n}) = \sum_{i\in \I_{0, n-1}} (bx_1+x_3)^{i}x_3^{n-1-i}.
\end{align*}
But, for $  n\geq 1 $, 
$
	(bx_1+x_3)^{n} = \sum_{i\in\I_{0,n}}\binom{n}{i}x_1^ix_3^{n-i}\prod_{j\in \I_{i}} \left( b+(j-1)\, \frac{b-p}{2}\right).$ Therefore,
\begin{align*}
\partial_3(x_3^{n}) &=  \sum_{\ell\in\I_{0, n-1}} \!\! \bigg( \sum_{i\in\I_{0, \ell}}\binom{\ell}{i}x_1^ix_3^{\ell-i}\prod_{j\in \I_{i}} \! \left( b+(j-1)\, \frac{b-p}{2}\right) \bigg) x_3^{n-1-\ell} \\
&=    \sum_{i\in\I_{0, n-1}}\binom{n}{i+1} x_1^i x_3^{n-1-i} \prod_{j\in \I_{i}} \left(  b+(j-1)\, \frac{b-p}{2}\right).
\end{align*}
\epf

 \begin{lem}\label{lem:1.2_1_det}
 	For $ n\geq 1 $, define the matrix $ M_n  $ of order $ n+1 $ given by
 	{\small \begin{align} \label{eqn:1.2_1_def M_n}
 		 \left( \begin{array}{ccccccc}
 	\beta_n	& P_{1}^1 & 0 & 0 & \cdots & 0& 0  \\
 	\beta_n\beta_{n-1}	& P_{2}^1\beta_{n-1}  & P_{2}^2 & 0 &  \cdots & 0 & 0 \\
 	\prod\limits_{\mathclap{j\in\I_{n-2, n}}} \! \beta_j 	& P_{3}^1\beta_{n-1}\beta_{n-2} & P_{3}^2\beta_{n-2}  & P_3^3 & \cdots & 0 & 0   \\
 	\vdots	 & \vdots & \vdots & \vdots &   \ddots & \vdots & \vdots \\
 	\prod\limits_{\mathclap{j\in\I_n}} \! \beta_j & P_{n}^1\prod\limits_{\mathclap{j\in\I_{n-1}}} \!  \beta_j & P_{n}^2\prod\limits_{\mathclap{j\in\I_{n-2}}} \!  \beta_j &P_{n}^3\prod\limits_{\mathclap{j\in\I_{n-3}}} \!  \beta_j &  \cdots & P_{n}^{n-1}\beta_{1} & P_{n}^n   \\
 	\gamma_{n}	& \gamma_{n-1} & \gamma_{n-2}  &\gamma_{n-3}  & \cdots & \gamma_{1} & 1  
 		\end{array}\right),
 	\end{align}}
where $ \beta_j $ and $ \gamma_j $ are as in \eqref{eqn:1.2_1_gamma_and_pi}, and $ P_i^j = i!/(i-j)! $.
Then,
\begin{align}\label{eqn:1.2_1_detM_n}
	\det M_n =\binom{n+1}{2} (p-b-2q) \prod\limits_{\mathclap{j\in\I_{n-1}}}  j!\, \beta_j.
\end{align}
 \end{lem}
\pf
After a cumbersome computation, we obtain that 
\begin{align*}
	\det M_n &= \bigg( \prod_{i\in \I_{n}} i!\bigg)  \sum_{\mathclap{j\in \I_{0, n}}}  \frac{(-1)^j\, \gamma_j}{(n-j)!}  \prod_{i\in \I_{j+1, n}} \!\!\!\!\! \beta_{i} , &&n\geq 1,
\end{align*} 
where $ \gamma_0 = 1 $.
Then \eqref{eqn:1.2_1_detM_n} follows since
\begin{align*}
  \sum_{j\in \I_{n-\ell, n}} \!\!\! \frac{(-1)^j\, \gamma_j}{(n-j)!}  \prod_{i\in \I_{j+1, n}} \!\!\!\!\! \beta_{i} &= \frac{n\, (-1)^{n-\ell}}{\ell! \, (n-\ell)}\gamma_{n-\ell}\prod_{\mathclap{i\in \I_{n-\ell, n-1}}} \beta_{i}, &&\ell\in \I_{0, n-1}.
\end{align*} 
\epf

\begin{pro} \label{prop:case1,2_azul}
	Assume that $ t= 1 $ and $ b\neq  p-2q$. 
	
	\label{case1:prop1.2_1} \rm{(a)} If $ \beta_n \neq 0 
	 $ for all $ n\in \mathbb{N} $, then $ \J(V) = \langle\eqref{eqn:1.2_1_relation degree 2}, \,\eqref{eqn:1.2_1_relation zn zn-1} \rangle  $, $ B_\infty 
	 $ is a PBW-basis of $ \bB(V) $ and $ \GK \bB(V)= \infty $.
	 
	\label{case2:prop1.2_1} \rm{(b)} Otherwise, there is a unique $ N\in\mathbb{N} $ such that $ \beta_N = 0 $; then $ \J(V) = \langle\eqref{eqn:1.2_1_relation degree 2}, \,\eqref{eqn:1.2_1_relation zn zn-1}, \, z_{N+1} + (N+1)(\frac{b+p}{2}) \, x_1z_N  \rangle  $, $ B_{N+3}=\{x_1^{a_1}x_2^{a_2}z_1^{b_1}z_2^{b_2} \cdots z_N^{b_N}x_3^{a_3}:  a_i, b_j\geq 0 \} $ is a PBW-basis of $ \bB(V) $ and $ \GK \bB(V)= N+3 $.
\end{pro}
\pf Write $ \bB = T(V)/I $ where $ I = \langle\eqref{eqn:1.2_1_relation degree 2}, \,\eqref{eqn:1.2_1_relation zn zn-1} \rangle $. We already proved that $ I\subset\J(V) $ and $ B_\infty $ generates linearly $ \bB $. Let $ \pi:\toba \to \toba(V) $ be the natural projection.

\smallbreak
 
\noindent  \hyperref[{case1:prop1.2_1}]{\rm{(a)}} Assume that $ \beta_n \neq 0 $ for all $ n\in \mathbb{N} $. Suppose that $ \pi $ is not an isomorphism and pick a linear homogeneous relation of minimal degree $ d \geq 3 $
 \begin{align*}
 	0\neq r = \sum\lambda_{a_1,a_2,a_3}^{b_1,\cdots, b_n}\,  x_1^{a_1}x_2^{a_2}z_1^{b_1}z_2^{b_2} \cdots z_n^{b_n}x_3^{a_3}\in\Bbbk B_\infty.
 \end{align*}
 We also denote $ \lambda_{a_1,a_2,a_3}^{b_1,\cdots, b_n}$ by $ \lambda_y $,  if $ y = x_1^{a_1}x_2^{a_2}z_1^{b_1} \cdots z_n^{b_n}x_3^{a_3}. $
 
Order the monomials $ B_\infty $ via the colexicographic order $ \prec $ considering that
\begin{align*}
	x_1 \prec x_2 \prec z_1\prec z_2\prec \cdots \prec z_n \prec \cdots\prec x_3.
\end{align*}
See \S \ref{subsec:colex_order} for details and notation.
Let $\bX = x_1^{\tilde{a}_1}x_2^{\tilde{a}_2}z_1^{\tilde{b}_1}z_2^{\tilde{b}_2} \cdots z_n^{\tilde{b}_n}x_3^{\tilde{a}_3} $ be the maximal term of the relation $ r $.
In particular, $\lambda_{\bX} \neq 0$.

If $ \tilde{a}_1\neq 0, $ then the maximal term of $ \partial_1 (\bX) $ is 
\begin{align*}
	\max\partial_1 (\bX) =  x_1^{\tilde{a}_1-1}x_2^{\tilde{a}_2}z_1^{\tilde{b}_1}z_2^{\tilde{b}_2} \cdots z_n^{\tilde{b}_n}x_3^{\tilde{a}_3}.
\end{align*}
Define $ B_{\prec \bX}^d = \{y\in B_\infty: \deg(y) = d \text{ and } y\prec \bX  \}$. Observe that the terms of $ \partial_1(y) $ are lower than $ \max\partial_1 (\bX) $ for all $ y\in B_{\prec \bX}^d $. 
Thus, the term $ \max\partial_1 (\bX) $ appears just one time in $ \partial_1 (r) = 0 $ what leads to $ \lambda_{\bX} = 0 $ what contradicts the maximality of $ \bX $. 
Hence $ \tilde{a}_1= 0 $ and $\bX = x_2^{\tilde{a}_2}z_1^{\tilde{b}_1}z_2^{\tilde{b}_2} \cdots z_n^{\tilde{b}_n}x_3^{\tilde{a}_3} $.

Similarly, $ \tilde{a}_2= 0 $; otherwise, $ \max\partial_2 (\bX) =  x_2^{\tilde{a}_2-1}z_1^{\tilde{b}_1}z_2^{\tilde{b}_2} \cdots z_n^{\tilde{b}_n}x_3^{\tilde{a}_3} $ and we apply the previous argument again.
Consequently, $\bX = z_1^{\tilde{b}_1}z_2^{\tilde{b}_2} \cdots z_n^{\tilde{b}_n}x_3^{\tilde{a}_3} $.

Let $ B^{ < d} = \{y\in B_\infty: \deg(y) < d \}$.
By the minimality of $ d $,  $ B^{<d} $ is a basis of $ \Bbbk B^{<d} $; then write $ (B^{ < d})^* = \{y^*\}_{y\in B^{ < d}}  $ for the dual basis of $ B^{ < d} $.

Suppose $ \tilde{b}_1\neq 0 $. Then $ \max\partial_1 (\bX) =  x_2z_1^{\tilde{b}_1-1}z_2^{\tilde{b}_2} \cdots z_n^{\tilde{b}_n}x_3^{\tilde{a}_3} $. Differently from the situation above, there is a unique $ \omega_{1,1}= x_1x_2 z_1^{\tilde{b}_1-1}z_2^{\tilde{b}_2}\cdots z_n^{\tilde{b}_n}x_3^{\tilde{a}_3} \in B_{\prec \bX}^d$ such that one term of $ \partial_1 (\omega_{1,1}) $ is not lower than $ \max\partial_1 (\bX) $. Furthermore, we have that $  \max\partial_1 (\omega_{1,1})  =  \max\partial_1 (\bX) $, whence
\begin{align}\label{eqn:1.2_1_M_1_linha 1}
\begin{split}
	0&=(\max\partial_1 (\bX))^* (\partial_1(r))\\
	&=(\max\partial_1 (\bX))^* (\lambda_{\bX}\partial_1(\bX)+\lambda_{\omega_{1,1}}\partial_1(\omega_{1,1})) = \tilde{b}_1\beta_1\, \lambda_{\bX}+\lambda_{\omega_{1,1}}.
\end{split}
\end{align}

On the other hand, $ \max\partial_2 (\omega_{1,1}) = x_1z_1^{\tilde{b}_1-1}z_2^{\tilde{b}_2}\cdots z_n^{\tilde{b}_n}x_3^{\tilde{a}_3} $ and, for all $ y\in B_{\prec \bX}^d - \{\omega_{1,1}\} $, the terms of $ \partial_2(y) $ are lower than $ \max\partial_2 (\omega_{1,1}) $. In particular,
\begin{align}\label{eqn:1.2_1_M_1_linha 2}
\begin{split}
0&=(\max\partial_2 (\omega_{1,1}))^* (\partial_2(r))\\
&=(\max\partial_2 (\omega_{1,1}))^* (\lambda_{\bX}\partial_2(\bX)+\lambda_{\omega_{1,1}}\partial_2(\omega_{1,1})) = \tilde{b}_1\gamma_1\, \lambda_{\bX}+\lambda_{\omega_{1,1}}.
\end{split}
\end{align}

Equations \eqref{eqn:1.2_1_M_1_linha 1} and \eqref{eqn:1.2_1_M_1_linha 2} give rise to a homogeneous system whose associated matrix is 
	$ \left( \begin{array}{cc}
	\tilde{b}_1\beta_1	& 1 \\
	\tilde{b}_1\gamma_{1}	&  1  
	\end{array}\right). $
	By Lemma \ref{lem:1.2_1_det}, the determinant of such matrix is $ \tilde{b}_1(p-b-2q) \neq 0 $ what implies that $ \lambda_{\bX}=0 $, a contradiction. Hence $ \tilde{b}_1 = 0 $ and $\bX = z_2^{\tilde{b}_2} \cdots z_n^{\tilde{b}_n}x_3^{\tilde{a}_3} $.

	Inductively, assume that $ \tilde{b}_i = 0 $ for $ i\in \I_{j-1} $ and suppose  $ \tilde{b}_j \neq 0 $. 
	Write $ \omega_{j,i}= x_1^iz_{j-i}z_j^{\tilde{b}_j-1}z_{j+1}^{\tilde{b}_{j+1}} \cdots z_n^{\tilde{b}_n}x_3^{\tilde{a}_3}, $ for $ i\in \I_{0, j} $.
Note that $ \omega_{j,0} = \bX $ and $ \omega_{j,i}\in B_{\prec \bX}^d, i\in \I_j $. 
It holds that $ \max\partial_1^\ell (\bX) =  z_{j-\ell}z_j^{\tilde{b}_j-1}z_{j+1}^{\tilde{b}_{j+1}} \cdots z_n^{\tilde{b}_n}x_3^{\tilde{a}_3}, \ell\in\I_j $.
		
	Moreover, if $ i> \ell $, then all the terms of $ \partial_1^\ell (\omega_{j,i}) $ are lower than $  \max\partial_1^\ell (\bX) $; if $ i\leq \ell $, then $ \max\partial_1^\ell (\omega_{j,i}) = \max\partial_1^\ell (\bX) $ and the list $ \{\omega_{j,i}\}_{i\in\I_{j}}  $ is exhaustive, that is, when we apply $ \partial_1^\ell $, $ \ell\in\I_j $, just the elements $ \{\omega_{j,i}\}_{i\in\I_{\ell}} $ of $ B_{\prec \bX}^d $ ``contribute'' in the direction of $ \max\partial_1^\ell (\bX) $.
	In particular, for all $ \ell\in\I_j $,
	\begin{align}\label{eqn:1.2_1_M_j_linha ell}
	\begin{split}
	0&=(\max\partial_1^\ell (\bX))^* (\partial_1^\ell(r))=(\max\partial_1^\ell (\bX))^* \bigg(\sum_{{i\in\I_{0, \ell}}}\lambda_{\omega_{j,i}}\partial_1^\ell(\omega_{j,i})\bigg) \\
	&= \lambda_{\bX} \, \tilde{b}_j \prod_{\mathclap{\theta\in \I_{j-\ell+1, j}}} \beta_\theta  + \sum_{\mathclap{{i\in\I_{\ell-1}}}}\lambda_{\omega_{j,i}} \, P_\ell^i   \prod_{\mathclap{\theta\in \I_{j-\ell+1, j-i}}} \beta_\theta. 
	\end{split}
	\end{align}

Further, $ \max\partial_2 (\omega_{j,j}) = x_1^jz_j^{\tilde{b}_j-1}z_{j+1}^{\tilde{b}_{j+1}}\cdots z_n^{\tilde{b}_n}x_3^{\tilde{a}_3} $ and the terms of $ \partial_2(y) $ are lower than $ \max\partial_2 (\omega_{j,j}) $ for all $ y\in B_{\prec \bX}^d - \{\omega_{j,i}\}_{i\in\I_j} $. Then, 
\begin{align}\label{eqn:1.2_1_M_j_linha j+1}
\begin{split}
0&=(\max\partial_2 (\omega_{j,j}))^* (\partial_2(r))=(\max\partial_2 (\omega_{j,j}))^* \bigg(\sum_{{i\in\I_{0, j}}}\lambda_{\omega_{j,i}}\partial_2(\omega_{j,i})\bigg) \\
&= \tilde{b}_j\gamma_j\, \lambda_{\bX}+\sum_{\mathclap{i\in\I_{j-1}}}\gamma_{j-i} \lambda_{\omega_{j,i}} + \lambda_{\omega_{j,j}}.
\end{split}
\end{align}

The matrix associated to the homogeneous system obtained from equations \eqref{eqn:1.2_1_M_j_linha ell} and \eqref{eqn:1.2_1_M_j_linha j+1} is equals to the matrix $ M_j $ \eqref{eqn:1.2_1_def M_n} with exception of the first column, which is the original one multiplied by the scalar $\tilde{b}_j$.
By Lemma \ref{lem:1.2_1_det}, the determinant of this matrix is $\tilde{b}_j \binom{j+1}{2} (p-b-2q) \prod_{{i\in\I_{j-1}}} \! i!\, \beta_i \neq 0 $, what implies that $ \lambda_{\bX}=0 $, a contradiction. 
Hence $ \tilde{b}_j = 0 $.

After the whole inductive process, we get
$\bX = x_3^{\tilde{a}_3} $.
As $ \tilde{a}_3= d \geq 3 $, then $ \max\partial_3 (\bX) = x_3^{\tilde{a}_3-1} $ by Lemma \ref{lem:1.2_1_partial 3}. 
Arguing as above, we have another con-tradiction.
Therefore, there are no more relations in $ \bB $ and $ \pi $ is, in fact, an isomorphism. In particular, $ B_\infty $ is a basis of $ \toba(V) $ and $ \GK\toba(V) =\infty $.

\smallbreak

\noindent  \hyperref[{case2:prop1.2_1}]{\rm{(b)}} Assume that $ \beta_N = 0 $ for some $ N\in \mathbb{N} $. 
An easy calculation shows that such $ N $ is unique.
Clearly $ \partial_i(z_{N+1} +(N+1)(\frac{b+p}{2})  \, x_1z_N) = 0 $ for all $ i\in\I_3 $.
Then consider the pre-Nichols algebra $ \widetilde{\bB} = T(V)/\widetilde{I} $ where $ \widetilde{I} = \langle\eqref{eqn:1.2_1_relation degree 2}, \break \eqref{eqn:1.2_1_relation zn zn-1},  \, z_{N+1} +(N+1)(\frac{b+p}{2})  \, x_1z_N \rangle $.
By induction on $ i $, we see that
\begin{align*}
	z_{N+i} &= \frac{(N+i)!}{N!} \left(\frac{-b-p}{2}\right)^i  x_1^iz_N \text{ holds in } \widetilde{\bB}, && i\geq 1.
\end{align*}
In particular, $ B_{N+3}=\{x_1^{a_1}x_2^{a_2}z_1^{b_1}z_2^{b_2} \cdots z_N^{b_N}x_3^{a_3}:  a_i, b_j\geq 0 \}  $ generates linearly $ \widetilde{\bB} $.

To get that $ \toba(V) = \widetilde{\bB} $, we proceed as in the proof of case \hyperref[{case1:prop1.2_1}]{\rm{(a)}}: we suppose the existence of another relation what leads to a contradiction. The fact that $ \beta_N = 0 $ does not affect the previous arguments because to verify that $ \tilde{b}_j = 0 $ for $ j\in \I_N $, we just need  that $ \beta_i \neq 0 $ for $ i\in \I_{N-1} $. 
\epf
%
%
%
%
%
%

\subsubsection{Case \ref{subcase:1,2_2}} \label{subsubsec:hit12_vermelho}
Assume that
$ t=-1 $ and $ b\neq  -p$.
We adapt the strategy used in \S \ref{subsubsec:hit12_azul} to classify the Nichols algebras under these conditions.

\smallbreak
By \eqref{eq:J2}, the quadratic relations of $ \bB(V) $ are
\begin{align}\label{eqn:1.2_2_relation degree 2}
&x_2x_1 + x_1x_2,  &&x_i^2,\quad i\in\I_2.
\end{align}
Set $ x_{31}=x_3x_1+x_1x_3 $ and note that $ x_{31}x_1 = x_{31}x_1 $.
Further, the relations
\begin{align}\label{eqn:1.2_2_relation degree 3}
& x_{31}x_2 -x_2x_{31}  &&\text{ and } &&x_3x_{31}-x_{31}x_3+(p-b)\, x_1x_{31}
\end{align}
hold in $ \toba(V) $.

For $ n\geq 1 $, define recursively $ z_n= \, (x_3-nb\, x_1)z_{n-1} +$
\begin{align*}
+\begin{cases}
 +z_{n-1}x_3 
+\left( \frac{n-1}{2}\right)!(p-q)(-p-b)^{ \frac{n-1}{2}} x_1x_2x_{31}^{ \frac{n-1}{2}}  &\text{if $ n $ is odd}, \\
 -z_{n-1}x_3 
-\left( \frac{n-2}{2}\right)!(p-q)(-p-b)^{ \frac{n-2}{2}} x_2x_{31}^{\frac{n}{2}} &\text{if $ n $ is even},
\end{cases}
\end{align*}
assuming that $ z_0:=x_2 $.
Then, the following properties
\begin{align}\label{eqn:1.2_2_properties1}
\begin{split}
z_nx_i =(-1)^{n+1} x_iz_n, 
\qquad\qquad\quad \,\,\,\,\,\,\,\,\,\,\,\, z_nx_{31} = x_{31}z_n, \!\!\!\!\!\!\!\!\!\!\!\!\qquad\qquad\qquad\\ 
z_nz_m =(-1)^{(n+1)(m+1)} z_nz_m,\quad\qquad\qquad z_{2n}^2 = 0 \qquad\qquad\qquad\\ 
c(z_n\otimes x_i) = (-1)^{n+1} x_i\otimes z_n, \qquad\qquad\qquad\qquad\qquad\quad \,\,\,\,\, \\ 
c(z_n\otimes x_3) = (-1)^{n+1}(x_3+(np+q)x_1)\otimes z_n, \qquad\qquad\quad 
\\ 
\partial_3(z_n) = 0,\qquad \quad
\partial_2(z_n) = \widetilde{\gamma}_n\, x_1x_{31}^{\frac{n-1}{2}}, \qquad\quad
\partial_1(z_n) = \widetilde{\beta}_n\, z_{n-1}, 
\end{split}
\end{align}
hold in $ \bB(V) $ for all $ n\geq 1 $, $ m\in\I_{n-1} $ and $ i\in \I_2 $, where
\begin{align}\label{eqn:1.2_2_gamma_and_pi}
\begin{split}
\widetilde{\gamma}_n&= \chi_{_{o}}(n)\left( \frac{n-1}{2}\right)!(q-p)\left( -p-b\right)^{ \frac{n-1}{2}} \text{ and } \,  \\\widetilde{\beta}_n&=\begin{cases}
-\left( \frac{n+1}{2}\, b+\frac{n-1}{2}\, p+q\right)&\text{if $ n $ is odd}, \\
-\frac{n}{2}\left(b+p\right)&\text{if $ n $ is even}.
\end{cases}
\end{split}
\end{align}

We observe that the relations
\begin{align}\label{eqn:1.2_2_relation zn zn-1_and_z_2n^2}
&z_{2n}z_{2n-1} - z_{2n-1}z_{2n}, &&z_{2n}^2, &&n\geq 1,
\end{align}
are new and proved by derivations.
We present the following identities that are necessary to show the properties \eqref{eqn:1.2_2_properties1}:
\begin{align*} 
\begin{split}
c(x_{31}\otimes x_i) &= x_i\otimes x_{31}, \\
\partial_j(x_{31}^n) &= -\delta_{j, 1} n(b+p)\, x_1x_{31}^{n-1}, \\
c(x_{31}^n\otimes x_3) &= (x_3+2np\, x_1)\otimes x_{31}^n, 
\end{split}
&&n\geq 0, \, i\in\I_2, \, j\in\I_3.
\end{align*}

\smallbreak

Consider the ideal $ I\subseteq\J(V) $ generated by \eqref{eqn:1.2_2_relation degree 2}, \eqref{eqn:1.2_2_relation degree 3} and \eqref{eqn:1.2_2_relation zn zn-1_and_z_2n^2}, and set $ \bB = T(V)/I$. 
Thus, for $ m, i,j\geq 0, \, i$ odd and $ j $ even,
\begin{align*}
\begin{split}
	x_3x_{1} &= x_{31} - x_{1}x_3, \\
	x_3x_{31}^m &= x_{31}^mx_3 - m(p-b)\, x_1x_{31}^m, \\
x_3z_{i}^m &=  \, z_{i}^mx_3 + m\, z_{i}^{m-1}z_{i+1}+ m(i+1)b \, x_1z_i^m \\
&+m\, \left( \frac{i-1}{2}\right)!(p-q)\left( -p-b\right)^{\frac{i-1}{2}}x_2x_{31}^{\frac{i+1}{2}}z_i^{m-1}, \\
x_3z_{j} &= z_{j+1}-  z_{j}x_3 +  (j+1)b \, x_1z_j \\
&- \left( \frac{j}{2}\right)!(p-q)\left( -p-b\right)^{\frac{j}{2}}x_1x_2x_{31}^{\frac{j}{2}},
\end{split}	
\end{align*}
hold in $ \bB $.
In particular, it follows that 
\begin{align}\label{eqn:1.2_2_Binfinito} 
\widetilde{B}_\infty = \{x_1^{a_1}x_2^{a_2}x_{31}^{a_3}z_1^{b_1} \cdots z_n^{b_n}x_3^{a_4}: n\geq 1; \, a_i, b_{j}\geq 0;\, a_1, a_2, b_{2j} < 2\}
\end{align}
is a system of linear generators of $ \bB $.
From now on, we use the notation $ a_i, b_j $ admitting always that they are suitable, that is, as in \eqref{eqn:1.2_2_Binfinito}.

Observe that
\begin{align*}
\partial_i(z_nz_m)&= \partial_i(z_n) z_m + (-1)^{n+1} z_n\partial_i(z_m), \qquad  i\in\I_3, \, n, m\geq 1. 
\end{align*}
Then, $ \partial_i(z_n^{b_n})=  b_n \,\partial_i(z_n) z_n^{b_n-1}, $ for all $  i\in\I_3, n\geq 1. $
Hence, we get that 
\begin{align*}
\begin{split}
\partial_i(x_1^{a_1}&x_2^{a_2}x_{31}^{a_3}z_1^{b_1}\cdots   z_n^{b_n}) = 
\partial_i(x_1^{a_1})x_2^{a_2}x_{31}^{a_3}z_1^{b_1} \cdots z_n^{b_n}\\
&+(-1)^{a_1} x_1^{a_1}\partial_i(x_2^{a_2})x_{31}^{a_3}z_1^{b_1} \cdots z_n^{b_n} +(-1)^{a_1+a_2} x_1^{a_1}x_2^{a_2}\partial_i(x_{31}^{a_3})z_1^{b_1} \cdots z_n^{b_n} \\
&+ \cdots  + (-1)^{a_1+a_2+2a_3+\, \sum\limits_{\mathclap{i\in\I_{n-1}}} \, (i+1)\, b_{i}} \! \!  x_1^{a_1}x_2^{a_2}x_{31}^{a_3}z_1^{b_1} \cdots \partial_i(z_n^{b_n}),
\end{split}
\end{align*} 
for all $ i\in\I_3 $. 

\begin{lem}\label{lem:1.2_2_partial 3} For $ a_i, b_j $ as in \eqref{eqn:1.2_2_Binfinito} and $ n>0 $, the derivations  
	\begin{align*}
	&\partial_i(\bY x_3^{a_4}) = \partial_i(\bY)\, x_3^{a_4} + (-1)^\eta\,\bY(\partial_i(x_3^{a_4}) +\delta_{i, 1}\alpha\, \partial_3(x_3^{a_4}) ), \\
	&\partial_3(x_3^{n}) =\begin{cases}
	-\sum\limits_{\mathclap{i\in\I_{0,\frac{n-2}{2}}}}\binom{\frac{n}{2}}{i+1}\, x_1x_{31}^ix_3^{n-2i-2}\prod\limits_{\mathclap{j\in \I_{i+1}}} jb-(j-1)p & \text{if $ n $ is even}, \\
	 \ \ \,  \sum\limits_{\mathclap{i\in\I_{0,\frac{n-1}{2}}}}\binom{\frac{n-1}{2}}{i}\, x_{31}^ix_3^{n-2i-1}\prod\limits_{\mathclap{j\in \I_{i}}} jb-(j-1)p &  \\
	-\sum\limits_{\mathclap{i\in\I_{0, \frac{n-3}{2}}}}\binom{\frac{n-1}{2}}{i+1}\, x_1x_{31}^ix_3^{n-2i-2}\prod\limits_{\mathclap{j\in \I_{i+1}}} jb-(j-1)p & \text{if $ n $ is odd},\end{cases}
	\end{align*} 
hold in $ \bB $, where $ \bY= x_1^{a_1}x_2^{a_2}x_{31}^{a_3}z_1^{b_1} \cdots z_n^{b_n} $, $ \eta = a_1+a_2+2a_3+\sum_{i\in\I_n} (i+1)\, b_{i} $ and  $ \alpha = a_1p+a_2q+2a_3p+\sum_{j\in \I_n} b_j(jp+q) $.
\end{lem} 

\pf 
Similar to Lemma \ref{lem:1.2_1_partial 3}.
We just mention that, for $ n>0 $,
%
%
%
\begin{align*}
(-bx_1-x_3)^{n} &=
\begin{cases}
x_3^n +\sum\limits_{\mathclap{i\in\I_{\frac{n}{2}}}}\binom{\frac{n}{2}}{i}\, x_{31}^ix_3^{n-2i}\prod\limits_{\mathclap{j\in \I_{i}}} jb-(j-1)p & \text{if $ n $ is even}, \\
-x_3^n -\sum\limits_{\mathclap{i\in\I_{\frac{n-1}{2}}}}\binom{\frac{n-1}{2}}{i}\, x_{31}^ix_3^{n-2i}\prod\limits_{\mathclap{j\in \I_{i}}} jb-(j-1)p &  \\
-\sum\limits_{\mathclap{i\in\I_{0, \frac{n-1}{2}}}}\binom{\frac{n-1}{2}}{i}\, x_1x_{31}^ix_3^{n-2i-1}\prod\limits_{\mathclap{j\in \I_{i+1}}} jb-(j-1)p & \text{if $ n $ is odd}.\end{cases}  
\end{align*}
\epf

\begin{lem}\label{lem:1.2_2_det}
For $ n\geq 1 $, let $ \widetilde{M}_n = (a_{i j})_{i,\, j\, \in\I_{n+1}} $ be the matrix given by
\begin{align*}
a_{ij} = \begin{cases}
0 &\text{ if } j \geq i+2, \\
\prod\limits_{\mathclap{s\in\I_{\ell}}} \widetilde{\beta}_{2s} &\text{ if } j = i+1, \\
0 &\text{ if } i, j\in\I_n,\,  j\leq i \text{ and $ i, j $ are even,}\\
\binom{\ell}{\theta}  \prod\limits_{\mathclap{s\in\I_{\theta}}}  \widetilde{\beta}_{2s} \ \   \prod\limits_{\mathclap{s\in\I_{ i-j+1}}} \! \widetilde{\beta}_{n+2-j-s} &\text{ if } i, j\in\I_n,\,  j\leq i \text{ and $ i $ or $ j $ is odd,}\\
\chi_{_{o}}(j)\, \widetilde{\gamma}_{n+1-j}  &\text{ if } i=n+1 \text{ and } j\in\I_n,\\
(-1)^n &\text{ if } i= j=n+1.
\end{cases}
\end{align*}
where
$ \ell= \left\lfloor\frac{i}{2}\right\rfloor $, $ \theta= \left\lfloor\frac{j-1}{2}\right\rfloor $, $ \widetilde{\beta}_j $ and $ \widetilde{\gamma}_j $ are as in \eqref{eqn:1.2_2_gamma_and_pi}.
Then, for $ m\geq 1 $,
\begin{align*}
\det \widetilde{M}_{2m} &=(-1)^m \prod\limits_{\mathclap{j\in\I_{m}}}  \widetilde{\beta}_{2j-1}   \widetilde{\beta}_{2j}^{ \ 2(m-j)+1}, \\ 
\det \widetilde{M}_{2m-1} &=(-1)^{m} m \widetilde{\beta}_2 \prod\limits_{\mathclap{j\in\I_{m-1}}}  \widetilde{\beta}_{2j-1} \widetilde{\beta}_{2j}^{ \ 2(m-j)}.
\end{align*}
\end{lem}

\pf
Analogous to Lemma \ref{lem:1.2_1_det}.
\epf

\begin{pro} \label{prop:case1,2_vermelho}
	Assume that $ t= -1 $ and $ b\neq  -p$. 
	
	\label{case1:prop1.2_2} \rm{(a)} If $ \widetilde{\beta}_n \neq 0 
	$ for all $ n\in \mathbb{N} $, then $ \J(V) = \langle\eqref{eqn:1.2_2_relation degree 2}, \eqref{eqn:1.2_2_relation degree 3}, \eqref{eqn:1.2_2_relation zn zn-1_and_z_2n^2} \rangle  $, $ \widetilde{B}_\infty 
	 $ is a PBW-basis of $ \bB(V) $ and $ \GK \bB(V)= \infty $.
	
	\label{case2:prop1.2_2} \rm{(b)} Otherwise, there is a unique odd number $ N\in\mathbb{N} $ such that $ \widetilde{\beta}_N = 0 $; then $ \J(V) = \langle\eqref{eqn:1.2_2_relation degree 2}, \eqref{eqn:1.2_2_relation degree 3}, \eqref{eqn:1.2_2_relation zn zn-1_and_z_2n^2}, \,  z_{N+1} + (N+1)(\frac{b+p}{2}) \, x_1z_N  \rangle  $, $ \widetilde{B}_{\frac{N+5}{2}}=\{ x_1^{a_1}x_2^{a_2}x_{31}^{a_3}z_1^{b_1}z_2^{b_2} \cdots z_N^{b_N}x_3^{a_4}:  a_i, b_j\geq 0; \, a_1, a_2, b_{2j} < 2 \} $ is a PBW-basis of $ \bB(V) $ and  $ \GK \bB(V)= \frac{N+5}{2} $.
\end{pro}
\pf We follow the same strategy and notation adopted in Proposition \ref{prop:case1,2_azul}.
Set 
$ x_1 \prec x_2 \prec x_{31} \prec z_1\prec z_2 \prec \cdots\prec x_3 $
and $\bX = x_1^{\hat{a}_1}x_2^{\hat{a}_2}x_{31}^{\hat{a}_3}z_1^{\hat{b}_1}z_2^{\hat{b}_2}   \cdots z_n^{\hat{b}_n}x_3^{\hat{a}_4} $. 

By the same reason that $ \tilde{a}_1 $ vanishes, we have that $ \hat{a}_i = 0,\, i\in\I_3. $
We also get $ \hat{b}_j = 0 $ applying Lemma \ref{lem:1.2_2_det}; here, 
\begin{align*}
	\omega_{j,i}&= x_1^{i-2\left\lfloor\frac{i}{2}\right\rfloor}x_{31}^{\left\lfloor\frac{i}{2}\right\rfloor} z_{j-i}z_j^{\hat{b}_j-1}z_{j+1}^{\hat{b}_{j+1}} \cdots z_n^{\hat{b}_n}x_3^{\hat{a}_4}, &&  i\in\I_{0, j}.
\end{align*}
Hence, $\bX = x_3^{\hat{a}_4} $.
If $ b=0 $ and $ \hat{a}_4 $ is even, then $ \partial_3 (x_3^{\hat{a}_4}) = 0 $ by Lemma \ref{lem:1.2_2_partial 3}.
To avoid this problem, we see that $ (x_3^{\hat{a}_4-1})^*(\partial_1(x_3^{\hat{a}_4})) = (-1)^{\hat{a}_4+1} p  \left\lfloor\frac{\hat{a}_4}{2}\right\rfloor  $ since
\begin{align*}
\begin{split}
\partial_1(x_3^{n+1}) = & -(x_3+ b\, x_1)\partial_1(x_3^{n}) + (q-p) \, x_2\partial_2(x_3^{n}) \\
&-(p\, x_3+ k\, x_2+ a\, x_1)\partial_3(x_3^{n}),
\end{split}
\qquad\qquad n >0. 
\end{align*}

For part  \hyperref[{case2:prop1.2_2}]{\rm{(b)}}, use that the following identity holds in $ \widetilde{\bB} $:
\begin{align*}
z_{N+i} &= \frac{\left\lfloor\frac{N+i}{2}\right\rfloor!}{\left(\frac{N-1}{2}\right)!} (-b-p)^{\left\lfloor\frac{i+1}{2}\right\rfloor}  x_1^{i-2\left\lfloor\frac{i}{2}\right\rfloor}x_{31}^{\left\lfloor\frac{i}{2}\right\rfloor}z_N, && i\geq 1.
\end{align*}

\epf

\subsubsection{Case \ref{subcase:1,2_3}} \label{subsubsec:hit12_remaining}

Here we classify the Nichols algebras associated to $ \hit_{1, 2} $ such that $ t^2=1 $ and are not covered by cases \ref{subcase:1,2_1} and \ref{subcase:1,2_2}.

\begin{pro} \label{prop:case1,2_remaining}
	Assume that  $ t= 1 $ and $ b=  p-2q$, or $ t= -1 $ and $ b=  -p$.
	Then, the Nichols algebras are as in Table \ref{tab:h12_3}, where
\begin{align}\label{eqn:J_1,2_1}
&\lg  x_2 x_1 - x_1x_2, \, x_3x_1 - x_1x_3 + q \, x_1^2, \, x_3x_2 - x_2x_3 + (2q-p) \, x_1x_2  \rg;  \\  \label{eqn:J_1,2_2}
& \ \ \ \ \lg  x_2x_1 + x_1x_2, \, x_3x_1 + x_1x_3, \, x_3x_2 + x_2x_3 + p\, x_1x_2, \, x_1^2,\, x_2^2 \rg;\\ \label{eqn:J_1,2_3}
& \lg  x_2x_1 + x_1x_2, \, x_3x_1 + x_1x_3, \, x_3x_2 + x_2x_3 + p\, x_1x_2, \,  \\
\notag   & \ \ \ \ \ \ \ \ \ \ \ \ \ \ \ \ \ \ \ \ \ \ \ \ \ \ \ \ \ \ \ \ \ \ \ \ \ \ \ \ \ \ \,  x_1^2, \,x_2^2, \, x_3^2 +k\, x_1x_2 + p\,x_1x_3 \rg.
\end{align}

\begin{table}[ht]
	\caption{Nichols algebras of type $ \hit_{1, 2} \, \ref{subcase:1,2_3}$}\label{tab:h12_3}
	\begin{center}
		\begin{tabular}{| c | c | c | c | c | c | c |}\hline
			Case & $ t $ & $ b $ & $ a $ &  $\J(V)$ &  Basis &$\GK$
			\\\hline
			 \label{case1:12} \rm{(a)}   & $ 1 $ & $ p-2q $ &  &  \eqref{eqn:J_1,2_1} & $ B_3 $ \eqref{eqn:B_3} & $ 3 $ \\\hline
			\label{case2:12} \rm{(b)} & $ -1 $ & $ -p $ & $ \neq -p^2 $   & \eqref{eqn:J_1,2_2} & $ B_1 $ \eqref{eqn:B_1} & $ 1 $ \\\hline
			\label{case3:12} \rm{(c)} & $ -1 $ & $ -p $ & $ -p^2 $  & \eqref{eqn:J_1,2_3} & $ B_0 $ \eqref{eqn:B_0} & $ 0, \, \dim =8 $ \\\hline
		\end{tabular}
	\end{center}
\end{table}
\end{pro}

\pf By \eqref{eq:J2}, the quadratic relations hold. 
Now we proceed the analysis of each case enumerated above separately.

\smallbreak

\noindent Case \hyperref[case1:12]{\rm{(a)}}: 
This case follows analogously to the case $ t=1 $ in Proposition \ref{prop:case1,1}.
For completeness, we present here just some of the necessary identities: 
\begin{align*}
	x_3x_1^{a_1} = x_1^{a_1}&x_3 - a_1q \, x_1^{a_1+1},\qquad\qquad  x_3x_2^{a_2} = x_2^{a_2}x_3 + a_2(p-2q) \, x_1x_2^{a_2},  \\
	&  c(x_1^{a_1}x_2^{a_2}\ot x_3) = (x_3 +  (a_1p+a_2q )x_1)\ot x_1^{a_1}x_2^{a_2}.
\end{align*}

\smallbreak

\noindent Case \hyperref[case2:12]{\rm{(b)}}:  
Write $ \bB = T(V)/I $ for the pre-Nichols algebra where $ I $ is the ideal \eqref{eqn:J_1,2_2}. 
It is clear that $B_1 =\{x_1^{a_1}x_2^{a_2}x_3^{a_3}: \, a_i\geq 0; \,  a_1, a_2 < 2 \}$ generates linearly $ \bB $.
Observe that the following derivations hold in $ \bB $:
\begin{align*}
\partial_1(x_1^{a_1}x_2^{a_2}x_3^{a_3}) &= \delta_{a_1, 1}\, x_2^{a_2}x_3^{a_3} -\delta_{a_2, 0}(-1)^{a_3+a_1}\left\lfloor (a_3+a_1)/2 \right\rfloor  p \, x_1^{a_1}x_3^{a_3-1} \\ 
&+ \delta_{a_2, 1}(-1)^{a_3} (\left\lfloor (a_3+a_1)/2 \right\rfloor p+\chi_{_{o}}(a_3) q)\, (-x_1)^{a_1}x_2x_3^{a_3-1} \\
&  -\left\lfloor a_3/2 \right\rfloor k \, ((-x_1)^{a_1}x_2^{a_2+1}x_3^{a_3-2}+\chi_{_{o}}(a_3) q  \, x_1^{a_1+1}x_2^{a_2+1}x_3^{a_3-3})   \\
&  +\left\lfloor a_3/2 \right\rfloor (\left\lfloor (a_3-2)/2 \right\rfloor p^2+a_2 pq -  a)  x_1^{a_1+1}x_2^{a_2}x_3^{a_3-2}, \\
\partial_2(x_1^{a_1}x_2^{a_2}x_3^{a_3})&= \ \delta_{a_2, 1} (-x_1)^{a_1}x_3^{a_3} + \left\lfloor a_3/2\right\rfloor k\, x_1^{a_1+1}x_2^{a_2}x_3^{a_3-2}, \\
\partial_3(x_1^{a_1}x_2^{a_2}x_3^{a_3})&= \ (-1)^{a_1+a_2}\chi_{_{o}}(a_3)\, x_1^{a_1}x_2^{a_2}x_3^{a_3-1}+ \left\lfloor a_3/2 \right\rfloor p\, x_1^{a_1+1}x_2^{a_2}x_3^{a_3-2},
\end{align*}
for all $ a_1, a_2\in \I_{0,1} $ and $ a_3\geq 0 $.

Suppose that the natural projection  $ \pi: \bB \to \toba (V) $ is not an isomorphism. 
Let 
$0\neq r = \lambda_{1}\, x_3^{n}+\lambda_{2}\, x_1x_3^{n-1}+ \lambda_{3}\, x_2x_3^{n-1}+ \lambda_{4}\, x_1x_2x_3^{n-2} \in \ker \pi$
a relation of minimal degree $ n\geq 3 $. 

If $ n $ is even, then we obtain  $ \lambda_{3} = 0 $ and $ 2\, \lambda_{4} = nk\, \lambda_{1} $ since $ \partial_2(r)=0 $. From $ \partial_3(r)=0 $, we get $ 2\, \lambda_{2} = np\, \lambda_{1} $. However $0 = 2\, \partial_1(r) = -n(p^2 + a)\lambda_{1}\, x_1x_3^{n-2}$ what implies $ \lambda_{1}=0 $ because $ a\neq -p^2 $.

If $ n $ is odd, then  $ r=0 $ using the following strategy: we seek a term that shows up just one time in some of the equations $ \partial_i (r)=0 $, $ i\in\I_3 $, and we use the linear independence given by the minimality of $ n $ to vanish the respective $ \lambda_j $. This argument is applied until we get all $ \lambda_j = 0$.

Hence $ r = 0 $ in both circumstances and, in particular, $ \pi $ is an isomorphism.
\smallbreak

\noindent Case \hyperref[case3:12]{\rm{(c)}}: Analogous to the case $ t=-1 $ of Proposition \ref{prop:case1,1}.
\epf

\subsection{Case $ \hit_{1, 3}$}\label{subsec:hit13}

Let $ c $ be the braiding associated to the solution of QYBE  
$$ \hit_{1, 3} = \left(\begin{array}{ccc|ccc|ccc}
t & \cdot & ta & \cdot & \cdot & \cdot & -ta & \cdot & -tab \\ 
\cdot & t & \cdot & \cdot & \cdot & \cdot & \cdot & -tb & -tp \\ 
\cdot & \cdot & t & \cdot & \cdot & \cdot & \cdot & \cdot & -tb  \\ \hline
\cdot & \cdot & \cdot & t & \cdot & tb & \cdot & \cdot & tp \\ 
\cdot & \cdot & \cdot & \cdot & t & \cdot & \cdot & \cdot & tq \\ 
\cdot & \cdot & \cdot & \cdot & \cdot & t & \cdot & \cdot & \cdot \\ \hline
\cdot & \cdot & \cdot & \cdot & \cdot & \cdot & t & \cdot & tb  \\ 
\cdot & \cdot & \cdot & \cdot & \cdot & \cdot & \cdot & t & \cdot  \\ 
\cdot & \cdot & \cdot & \cdot & \cdot & \cdot & \cdot & \cdot & t 
\end{array} \right). $$
Apply \eqref{eqn:QYBE_to_braiding} for the explicit presentation of $ c $.

%
%

\begin{pro} \label{prop:case1,3}
	If $ t^2 \neq 1 $, then there are no quadratic relations. Otherwise, the Nichols algebras are as in Table \ref{tab:h13}, where
	\begin{align}\label{eqn:J_1,3_1}
	\lg  x_2x_1 - x_1x_2, \, x_3x_1 - x_1x_3 - a \, x_1^2, &\, x_3x_2 - x_2x_3 -b \, x_1x_2  \rg; \\\label{eqn:J_1,3_2}
	\lg  x_2x_1 + x_1x_2, \, x_3x_1 + x_1x_3, \, x_3x_2 &+ x_2x_3 -b\, x_1x_2,\, x_1^2,\, x_2^2 \rg;
	 \\\label{eqn:J_1,3_3}
	 \lg  x_2x_1 + x_1x_2, \, x_3x_1 + x_1x_3, \, x_3x_2 + & x_2x_3 -b\, x_1x_2,\,   \\
	 \notag  &  x_1^2, \, x_2^2, \, x_3^2 -b\, x_1x_3 -p\, x_1x_2  \rg.
	\end{align}  
	
	\begin{table}[ht]
		\caption{Nichols algebras of type $ \hit_{1, 3} $}\label{tab:h13}
		\begin{center}
			\begin{tabular}{| c | c | c | c | c | c |}\hline
				 Case  & $ t $ & $ q $ & $\J(V)$ &  Basis &$\GK$
				\\\hline
				 \label{case1:13} \rm{(a)} & $1$ &   & \eqref{eqn:J_1,3_1} & $ B_3 $ \eqref{eqn:B_3} & $ 3 $ \\\hline
				\label{case2:13} \rm{(b)}   & $-1$ & $\neq 0$ & \eqref{eqn:J_1,3_2} & $ B_1 $ \eqref{eqn:B_1} & $ 1 $ \\\hline
				\label{case3:13} \rm{(c)}   &$-1$ & $ 0$ & \eqref{eqn:J_1,3_3} & $ B_0 $ \eqref{eqn:B_0} & $ 0, \, \dim = 8 $ \\\hline
			\end{tabular}
		\end{center}
	\end{table}	

\end{pro}

\pf We get the quadratic relations through  \eqref{eq:J2}. 
Next we study the cases individually.

\smallbreak

\noindent Case \hyperref[case1:13]{\rm{(a)}}:
 It follows similarly to case $ t=1 $ in Proposition \ref{prop:case1,1}.
 Some of the necessary identities to prove it are: 
\begin{align*}
x_3x_1^{a_1} = x_1^{a_1}x_3 + a_1a \, x_1^{a_1+1}, \qquad\qquad\quad x_3x_2^{a_2} &= x_2^{a_2}x_3 + a_2b \, x_1x_2^{a_2},  \\
c(x_1^{a_1}x_2^{a_2}\ot x_3) = (x_3 - (a_1a +a_2b )x_1) &\ot x_1^{a_1}x_2^{a_2}.
\end{align*}

\smallbreak

\noindent Case \hyperref[case2:13]{\rm{(b)}}:
Imitate case \hyperref[case2:12]{\rm{(b)}}  of Proposition \ref{prop:case1,2_remaining}.
For $ a_1, a_2\in \I_{0,1} $ and $ a_3\geq 0 $, the corresponding derivations are: 
\begin{align*}
\partial_1(x_1^{a_1}x_2^{a_2}x_3^{a_3}) =& \ \delta_{a_1, 1} x_2^{a_2}x_3^{a_3} + p \left\lfloor {a_3}/{2}\right\rfloor (-x_1)^{a_1}x_2^{a_2+1}x_3^{a_3-2} \\
& \hspace*{-9pt}  + (\left\lfloor {a_3}/{2}\right\rfloor (-1)^{a_3} b - \chi_{_{o}}(a_3) (a_1a+a_2b))\, (-x_1)^{a_1}(-x_2)^{a_2}x_3^{a_3-1} \\ 
&  \hspace*{-9pt}+ b \left\lfloor {a_3}/{2}\right\rfloor  (a\chi_{_{e}}(a_3) + b(a_2+\left\lfloor {(a_3-1)}/{2}\right\rfloor))\, x_1^{a_1+1}x_2^{a_2}x_3^{a_3-2} \\ 
& \hspace*{-9pt} -pa \left\lfloor {a_3}/{2}\right\rfloor \chi_{_{o}}(a_3) \, x_1^{a_1+1}x_2^{a_2+1}x_3^{a_3-3}, \\ 
\partial_2(x_1^{a_1}x_2^{a_2}x_3^{a_3})=& \ \delta_{a_2, 1} (-x_1)^{a_1}x_3^{a_3} - \left\lfloor {a_3}/{2} \right\rfloor q \, (-x_1)^{a_1}x_2^{a_2+1}x_3^{a_3-2}\\
& - \left\lfloor {a_3}/{2}\right\rfloor x_1^{a_1+1}x_2^{a_2}\left( p\, x_3 -  \chi_{_{o}}(a_3) qb \, x_2 \right) x_3^{a_3-3}, \\
\partial_3(x_1^{a_1}x_2^{a_2}x_3^{a_3})=& \ \chi_{_{o}}(a_3)\, (-x_1)^{a_1}(-x_2)^{a_2}x_3^{a_3-1} - \left\lfloor {a_3}/{2} \right\rfloor b \, x_1^{a_1+1}x_2^{a_2}x_3^{a_3-2}.
\end{align*}

\smallbreak

\noindent Case \hyperref[case3:13]{\rm{(c)}}: 
Analogous to case $ t=-1 $ of Proposition \ref{prop:case1,1}.
\epf

\subsection{Case $ \hit_{1, 4}$}\label{subsec:hit14}
Let $ c $ be the braiding associated to the solution of QYBE  
$$ \hit_{1, 4} = \left(\begin{array}{ccc|ccc|ccc}
t & \cdot & ta & \cdot & \cdot & \cdot & -ta & \cdot & -tab \\ 
\cdot & t & \cdot & \cdot & \cdot & \cdot & \cdot & t(a-2b) & tp \\ 
\cdot & \cdot & t & \cdot & \cdot & \cdot & \cdot & \cdot & -tb  \\ \hline
\cdot & \cdot & \cdot & t & \cdot & t(2b-a) & \cdot & \cdot & \cdot \\ 
\cdot & \cdot & \cdot & \cdot & t & \cdot & \cdot & \cdot & \cdot \\ 
\cdot & \cdot & \cdot & \cdot & \cdot & t & \cdot & \cdot & \cdot \\ \hline
\cdot & \cdot & \cdot & \cdot & \cdot & \cdot & t & \cdot & tb  \\ 
\cdot & \cdot & \cdot & \cdot & \cdot & \cdot & \cdot & t & \cdot  \\ 
\cdot & \cdot & \cdot & \cdot & \cdot & \cdot & \cdot & \cdot & t 
\end{array} \right). $$
See \eqref{eqn:QYBE_to_braiding} for details on the presentation of $ c $.


\begin{pro} \label{prop:case1,4}
	If $ t^2 \neq 1 $, then there are no quadratic relations. Otherwise, the Nichols algebras are as in Table \ref{tab:h14}, where 
\begin{align}\label{eqn:J_1,4_1}
\lg  x_2x_1 - x_1x_2, \, x_3x_1 - x_1x_3 &- a \, x_1^2, \, x_3x_2 - x_2x_3 + (a - 2b) \, x_1x_2  \rg; \\ \label{eqn:J_1,4_2}
\lg  x_2x_1 + x_1x_2, \, x_3x_1 + x_1x_3,& \, x_3x_2 + x_2x_3 +(a - 2b) \, x_1x_2,\, x_1^2,\, x_2^2 \rg; \\\label{eqn:J_1,4_3}
\lg  x_2x_1 + x_1x_2, \, x_3x_1 + x_1x_3, &\, x_3x_2 + x_2x_3 +(a - \,  \,2b) \, x_1x_2,\, \\
\notag & \qquad\qquad\qquad \ \ \ \ \  x_1^2,\, x_2^2, \, x_3^2 -b\, x_1x_3 \rg.
\end{align}
	
	\begin{table}[ht]
	\caption{Nichols algebras of type $ \hit_{1, 4} $}\label{tab:h14}
	\begin{center}
		\begin{tabular}{| c | c | c | c | c | c |}\hline
			Case & $ t $ & $ p $ & $\J(V)$ &  Basis &$\GK$
			\\\hline
			\label{case1:14} \rm{(a)} & $1$ &   & \eqref{eqn:J_1,4_1} & $ B_3 $ \eqref{eqn:B_3} & $ 3 $ \\\hline
			\label{case2:14} \rm{(b)} &$-1$ & $\neq 0$ & \eqref{eqn:J_1,4_2} & $ B_1 $ \eqref{eqn:B_1} & $ 1 $ \\\hline
			\label{case3:14} \rm{(c)} &$-1$ & $ 0$ & \eqref{eqn:J_1,4_3} & $ B_0 $ \eqref{eqn:B_0} & $ 0, \, \dim = 8 $ \\\hline
		\end{tabular}
	\end{center}
\end{table}		
	
\end{pro}

\pf Similar to Proposition \ref{prop:case1,3}. 
For completeness, we present next just some of the necessary identities. For case \hyperref[case1:14]{\rm{(a)}}:
\begin{align*}
x_3x_1^{a_1} &= x_1^{a_1}x_3 + a_1a \, x_1^{a_1+1}, & x_3&x_2^{a_2} = x_2^{a_2}x_3 - a_2(a-2b) \, x_1x_2^{a_2},  \\
&c(x_1^{a_1}x_2^{a_2}\ot x_3) =  (x_3 +  & \!\!\!\!\! ( a_2 ( a - 2 b ) & - a_1 a ) x_1 ) \ot x_1^{a_1}x_2^{a_2},
\end{align*}
and for case \hyperref[case2:14]{\rm{(b)}}:
\begin{align*}
\partial_1(x_1^{a_1}x_2^{a_2}x_3^{a_3}) =& \ \delta_{a_1, 1} x_2^{a_2}x_3^{a_3} - p \left\lfloor {a_3}/{2}\right\rfloor (-x_1)^{a_1}x_2^{a_2+1}x_3^{a_3-2} \\
& \hspace*{-42pt}  + (\left\lfloor {a_3}/{2}\right\rfloor (-1)^{a_3} b - \chi_{_{o}}(a_3) (a_1a+a_2(2b-a)))\, (-x_1)^{a_1}(-x_2)^{a_2}x_3^{a_3-1} \\ 
&  \hspace*{-42pt}+ b \left\lfloor {a_3}/{2}\right\rfloor  (a(\chi_{_{e}}(a_3)-a_2)+b(2a_2+\left\lfloor {(a_3-1)}/{2}\right\rfloor))\, x_1^{a_1+1}x_2^{a_2}x_3^{a_3-2} \\ 
& \hspace*{-42pt} +pb \left\lfloor {a_3}/{2}\right\rfloor \chi_{_{o}}(a_3) \, x_1^{a_1+1}x_2^{a_2+1}x_3^{a_3-3}, \\ 
\partial_2(x_1^{a_1}x_2^{a_2}x_3^{a_3})=& \ \delta_{a_2, 1} (-x_1)^{a_1}x_3^{a_3}, \\
\partial_3(x_1^{a_1}x_2^{a_2}x_3^{a_3})=& \ \chi_{_{o}}(a_3)\, (-x_1)^{a_1}(-x_2)^{a_2}x_3^{a_3-1} - \left\lfloor {a_3}/{2} \right\rfloor b \, x_1^{a_1+1}x_2^{a_2}x_3^{a_3-2}.
\end{align*}
\epf

\subsection{Case $ \hit_{1, 5}$}\label{subsec:hit15}
Let $ c $ be the braiding associated to the solution of QYBE  
$$ \hit_{1, 5} = \left(\begin{array}{ccc|ccc|ccc}
t & t\ell & \cdot & -t\ell & -t\ell^2 & t(\ell q-k) & tp & tk & ta \\ 
\cdot & t & \cdot & \cdot & -t\ell & \cdot & \cdot & tq & tb \\ 
\cdot & \cdot & t & \cdot & \cdot & -t\ell & \cdot & \cdot & tp  \\ \hline
\cdot & \cdot & \cdot & t & t\ell & \cdot & \cdot & t(p-q) & -tb \\ 
\cdot & \cdot & \cdot & \cdot & t & \cdot & \cdot & \cdot & \cdot \\ 
\cdot & \cdot & \cdot & \cdot & \cdot & t & \cdot & \cdot & \cdot \\ \hline
\cdot & \cdot & \cdot & \cdot & \cdot & \cdot & t & t\ell & \cdot  \\ 
\cdot & \cdot & \cdot & \cdot & \cdot & \cdot & \cdot & t & \cdot  \\ 
\cdot & \cdot & \cdot & \cdot & \cdot & \cdot & \cdot & \cdot & t 
\end{array} \right). $$
Use  \eqref{eqn:QYBE_to_braiding} for the explicit presentation of $ c $.

%
%
%

The Nichols algebra associated to $ c $ has a quadratic relation if and only if $ t^2=1 $.
In the next two results, we compute the Nichols algebras when $ t=1 $ and $ p=2q $, and $ t=-1 $ and $ p=0 $, respectively.

\begin{pro} \label{prop:case1,5_t=1_and_p=2q}
	Assume that $ t=1 $ and $ p=2q $. Then
	\begin{align*}
	\J(V) = \lg  x_2x_1 - x_1x_2 - \ell \, x_1^2,  & \ x_3x_1 - x_1x_3 +  q \, x_1^2, \\
	&  \   x_3x_2 - x_2x_3 + \ell \, x_1x_3  + q \, x_1x_2 + (k-q\ell) \, x_1^2    \rg,
	\end{align*}	
$ B_3 \ \eqref{eqn:B_3}	$ is a PBW-basis of $ \bB(V) $ and $ \GK \bB(V)= 3 $.
\end{pro}

\pf The relations above follow by \eqref{eq:J2}. 
Let $ \bB = T(V)/I $ be the pre-Nichols algebra   where I is the ideal generated by these quadratic relations.
Observe that, for $ a_1, a_2\geq 0 $, the relations
\begin{align*}
\notag x_2x_1^{a_1} = x_1^{a_1}x_2   + a_1\ell \, x_1^{a_1+1}, \qquad \ \   x_3x_1^{a_1} = x_1^{a_1}x_3 - a_1q \, x_1^{a_1+1},\hspace{1.5cm} \\ 
x_3x_2^{a_2} \!  =   x_2^{a_2}x_3 - a_2 \ell  x_1x_2^{a_2-1}x_3 -  a_2 q   x_1 x_2^{a_2}  + (q\ell-k)  \sum_{\mathclap{i\in \I_{a_2}}} \frac{a_2! \, \ell^{i-1}}{(a_2-i)!   i}  x_1^{i+1}x_2^{a_2-i} \! ,
\end{align*}
hold in $ \bB $. They guarantee that $ B_3 $ is a system of linear generators of $ \bB $.

As $ c(x_1^{a_1}x_2^{a_2}\ot x_3)\in \bB_1\otimes \bB_{a_1+a_2} $, we have  it is equal to $  \sum_{i\in\I_3} x_i\otimes h_{i}$ for some  $ h_{i} \in \bB_{a_1+a_2} $. We prove that $ h_3 = \sum_{{i\in \I_{0, a_2}}} \frac{a_2 ! \, \ell^{i}}{(a_2-i) !} x_1^{a_1 + i}x_2^{a_2-i} $ by induction. Then,
\begin{align}\label{eqn:derivation2_case1_1,5}
\partial_3(x_1^{a_1}x_2^{a_2}x_3^{a_3}) &=  \sum_{\mathclap{i\in \I_{0,a_2}}} \dfrac{a_2! \, a_3 \, \ell^i }{(a_2-i)!}  \, x_1^{a_1+i}x_2^{a_2-i}x_3^{a_3-1}, \quad \, a_1, a_2, a_3\geq 0,
\end{align}
since $ \partial_3(x_3^{a_3}) = a_3x_3^{a_3-1} $.

Suppose that the natural projection  $ \pi: \bB \to \toba (V) $ is not injective. Pick a linear homogeneous relation  of minimal degree $ n\geq 3 $
\begin{align*}
r =  \sum_{\mathclap{a_1+a_2+a_3=n,  \, a_i\geq 0}}  \lambda_{a_1, a_2, a_3}\, x_1^{a_1}x_2^{a_2}x_3^{a_3}.
\end{align*}
By \eqref{eqn:derivation2_case1_1,5}, we have that
\begin{align}\label{eqn:partial_sum_prop_1.5}
\begin{split}
0 &= \partial_3(r) =  \sum_{ \mathclap{a_1+a_2+a_3=n, \,  a_i\geq 0}}  \lambda_{a_1, a_2, a_3}\, a_3 \sum_{\mathclap{i\in \I_{0, a_2}}} \dfrac{a_2! \, \ell^i}{(a_2-i)!}  \, x_1^{a_1+i}x_2^{a_2-i}x_3^{a_3-1} \\
&=  \sum_{{a_2\in \I_{0, n-1}}} \sum_{{i\in \I_{0, a_2}}} \sum_{ \substack{a_1+a_3=n-a_2 \\ a_1\geq 0, \, a_3\geq 1}} \! \! \! \! \! \! \! \! \! \lambda_{a_1, a_2, a_3}\,\dfrac{a_2! \, a_3 \, \ell^i}{(a_2-i)!}  \, x_1^{a_1+i}x_2^{a_2-i}x_3^{a_3-1}.
\end{split}
\end{align}

Note that the term $ x_2^{n-1} $ appears just one time in \eqref{eqn:partial_sum_prop_1.5}. 
By the minimality of $ n $, we get $ \lambda_{0, n-1, 1}=0 $. 
Now, for each $ j\in\I_{0, 1} $, note that the term $ x_1^{j}x_2^{n-2} x_3^{1-j} $ shows up only one time in \eqref{eqn:partial_sum_prop_1.5} and then we obtain $ \lambda_{j, n-2, 2-j}=0 $.
Inductively on $ k $, we have that $ \lambda_{j, n-k, k-j}=0, \, j\in\I_{0, k-1} $, for all $ k\in \I_n $.
In particular, we can rewrite the relation $ r $  as $  \sum \lambda_{a_1, a_2, 0}\, x_1^{a_1}x_2^{a_2}$.

Applying that
\begin{align*} 
\partial_2(x_1^{a_1}x_2^{a_2}) &= \sum_{i\in \I_{a_2}} \dfrac{a_2!\, \ell^{i-1}}{(a_2-i)!\, i}  \,x_1^{a_1 + i -1}x_2^{a_2-i}, \quad \, a_1\geq 0, \, a_2\geq 1,
\end{align*}
and repeating the previous procedure, we obtain that  $r=0$, a contradiction.
Hence, there are no more relations and $ \pi $ is an isomorphism.
\epf

\begin{pro} \label{prop:case1,5_t=-1_and_p_neq_0}
	Assume that $ t=-1 $ and $ p\neq 0 $.
	If $ a \neq b\ell $, then $ \J(V)= $
	\begin{align*}
	\lg  x_2x_1  &+ x_1x_2, \, x_3x_1 + x_1x_3,\, x_3x_2 + x_2x_3 -\ell \, x_1x_3   + q \, x_1x_2,  x_1^2,\, x_2^2 - \ell \, x_1x_2 \,  \rg,
	\end{align*}
$ B_1 \ \eqref{eqn:B_1}	$ is a PBW-basis of $ \bB(V) $ and $ \GK \bB(V)= 1 $.
If $ a = b\ell $, then 
\begin{align*}
\J(V) = \lg  x_2x_1 &+ x_1x_2, \, x_3x_1 + x_1x_3, \, x_3x_2 \, + \,  x_2x_3  - \ell \, x_1x_3   + q \, x_1x_2, \\
	\notag \,  & \qquad\qquad\qquad\qquad\qquad\qquad \ \ \ \ \, x_1^2,\, x_2^2 - \ell \, x_1x_2, \, x_3^2 + b \, x_1x_2 \,  \rg,
	\end{align*}
$ B_0 \ \eqref{eqn:B_0}	$ is a PBW-basis of $ \bB(V) $ and $ \dim \bB(V)= 8 $.
\end{pro}

\pf
By \eqref{eq:J2}, we get the relations above. 
If $ a=b\ell $, then it is clear that no more relations exist.

If $ a\neq b\ell $, this case follows the same lines of case \hyperref[case2:12]{\rm{(b)}} of Proposition \ref{prop:case1,2_remaining}.
We present here only the derivations: 
\begin{align*}
\partial_1(x_3^{a_3}) &= - \left\lfloor a_3/2\right\rfloor \left( a \, x_1x_3 + b \, x_2x_3 +   \chi_{_{o}}(a_3) \,  qb \, x_1x_2 \right) x_3^{a_3-3}, \\
\partial_1(x_1x_3^{a_3}) &=   x_3^{a_3} + \left\lfloor a_3/2\right\rfloor b \, x_1x_2x_3^{a_3-2}, \\
\partial_1(x_2x_3^{a_3}) &=  \left\lfloor a_3/2\right\rfloor (b\ell - a) \, x_1x_2x_3^{a_3 - 2} - \chi_{_{o}}(a_3)  \left( k \, x_1 + q \, x_2 \right)  x_3^{a_3 - 1}, \\
\partial_1(x_1x_2x_3^{a_3}) &=  x_2x_3^{a_3} + \ell \, x_1x_3^{a_3} + \chi_{_{o}}(a_3)\, q \, x_1x_2x_3^{a_3 - 1}, \\
%
%
\partial_2(x_1^{a_1}x_2^{a_2}x_3^{a_3})&=  \delta_{a_2, 1} x_1^{a_1}((-1)^{a_1}x_3\! +\!  \chi_{_{o}}(a_3) q x_1 \!)x_3^{a_3-1} \!\! + \! \left\lfloor a_3/2\right\rfloor \! b x_1^{a_1+1}x_2^{a_2}x_3^{a_3-2}\!\! , \\
\partial_3(x_1^{a_1}x_2^{a_2}x_3^{a_3})&=  \chi_{_{o}}(a_3)\, (-x_1)^{a_1}(-\ell \, x_1 - x_2)^{a_2}x_3^{a_3-1}.
\end{align*}
\epf

\subsection{Case $ \hit_{1, 6}$}\label{subsec:hit16}
Let $ c $ be the braiding associated to the solution of QYBE  
$$ \hit_{1, 6} = \left(\begin{array}{ccc|ccc|ccc}
t & ta & tp & -ta & -tab & -t(2pa+k) & -tp & tk & -tpq \\ 
\cdot & t & \cdot & \cdot & -tb & \cdot & \cdot & -tp & \cdot \\ 
\cdot & \cdot & t & \cdot & \cdot & -ta & \cdot & \cdot & -tq  \\ \hline
\cdot & \cdot & \cdot & t & tb & tp & \cdot & \cdot & \cdot \\ 
\cdot & \cdot & \cdot & \cdot & t & \cdot & \cdot & \cdot & \cdot \\ 
\cdot & \cdot & \cdot & \cdot & \cdot & t & \cdot & \cdot & \cdot \\ \hline
\cdot & \cdot & \cdot & \cdot & \cdot & \cdot & t & ta & tq  \\ 
\cdot & \cdot & \cdot & \cdot & \cdot & \cdot & \cdot & t & \cdot  \\ 
\cdot & \cdot & \cdot & \cdot & \cdot & \cdot & \cdot & \cdot & t 
\end{array} \right). $$
For the presentation of $ c $ apply \eqref{eqn:QYBE_to_braiding}.


\begin{pro} \label{prop:case1,6}
If $ t^2 \neq 1 $, then there are no quadratic relations. Otherwise, the Nichols algebras are as in Table \ref{tab:h16}, where 
\begin{align}\label{eqn:Nichols1,6_t=-1}
\lg  x_2x_1 + x_1x_2, \, x_3x_1 + x_1x_3, \, x_3x_2 \, + \, &  x_2x_3 -  a \, x_1x_3 - p \, x_1x_2, \\
\notag  & \, x_1^2,\, x_2^2 - b \, x_1x_2, \, x_3^2 -q\, x_1x_3 \rg; \\ \label{eqn:Nichols1,6_t=1}
\lg  x_2x_1 - x_1x_2 - a \, x_1^2,  \, x_3x_1 - x_1x_3 - & p\, x_1^2, \\
\notag  \, x_3x_2 - x_2x_3 + & a \, x_1x_3 -p \, x_1x_2 + (k+pa) \, x_1^2 \rg.
\end{align}

\begin{table}[ht]
	\caption{Nichols algebras of type $ \hit_{1, 6} $}\label{tab:h16}
	\begin{center}
		\begin{tabular}{| c | c | c | c | c | c |}\hline
			$ t $  & $\J(V)$ &  Basis &$\GK$
			\\\hline
			$-1$ &  \eqref{eqn:Nichols1,6_t=-1} & $ B_0 $ \eqref{eqn:B_0} & $ 0, \, \dim = 8 $ \\\hline
			$1$ &  \eqref{eqn:Nichols1,6_t=1} & $ B_3  $ \eqref{eqn:B_3} & $ 3 $ \\\hline
		\end{tabular}
	\end{center}
\end{table}		
\end{pro}

\pf 
The relations above hold by \eqref{eq:J2}. 
Clearly, there are no more relations if $ t=-1 $.
If $ t=1 $, then this case is similar to Proposition \ref{prop:case1,5_t=1_and_p=2q}.
Next, we present just some of necessary identities: 
\begin{align*}
&  x_2x_1^{a_1} = x_1^{a_1}x_2 + a_1a \, x_1^{a_1+1}, \qquad\quad x_3x_1^{a_1} = x_1^{a_1}x_3 + a_1p \, x_1^{a_1+1}, \\ 
 x_3x_2^{a_2} =& \ x_2^{a_2}x_3 - a_2 a  x_1x_2^{a_2-1}x_3 +  a_2 p  x_1 x_2^{a_2} \!  - (pa+k) \! \sum_{\mathclap{i\in\I_{a_2}}} \frac{a_2! \, a^{i-1}}{(a_2-i)!  i} x_1^{i+1} \! x_2^{a_2-i}\! , \\
&\partial_2(x_1^{a_1}x_2^{a_2}x_3^{a_3}) = \sum_{\mathclap{i\in\I_{a_2}}} {a_2 \choose i} \ \prod_{\mathclap{j\in\I_{0, i-2}}} (ja+b) \, x_1^{a_1-1+i}x_2^{a_2-i}x_3^{a_3}, \\ 
& \ \ \ \,  \partial_3(x_1^{a_1}x_3^{a_3}) = \sum_{\mathclap{i\in\I_{a_3}}} {a_3 \choose i}  \ \prod_{\mathclap{j\in\I_{0, i-2}}} (jp+q) \, x_1^{a_1-1+i}x_3^{a_3-i}.  
\end{align*}
\epf

\subsection{Case $ \hit_{1, 7}$}\label{subsec:hit17}
Let $ c $ be the braiding associated to the solution of QYBE  
$$ \hit_{1, 7} = \left(\begin{array}{ccc|ccc|ccc}
t & tk & tp & -tk & -tk\ell & td & ta & t(k(a-q)-d) & tb \\ 
\cdot & t & \cdot & \cdot & -t\ell & t(p-q) & \cdot & ta & \cdot \\ 
\cdot & \cdot & t & \cdot & \cdot & -tk & \cdot & \cdot & ta  \\ \hline
\cdot & \cdot & \cdot & t & t\ell & tq & \cdot & \cdot & \cdot \\ 
\cdot & \cdot & \cdot & \cdot & t & \cdot & \cdot & \cdot & \cdot \\ 
\cdot & \cdot & \cdot & \cdot & \cdot & t & \cdot & \cdot & \cdot \\ \hline
\cdot & \cdot & \cdot & \cdot & \cdot & \cdot & t & tk & tp  \\ 
\cdot & \cdot & \cdot & \cdot & \cdot & \cdot & \cdot & t & \cdot  \\ 
\cdot & \cdot & \cdot & \cdot & \cdot & \cdot & \cdot & \cdot & t 
\end{array} \right). $$
For details on the explicit presentation of $ c $ see \eqref{eqn:QYBE_to_braiding}.


The Nichols algebra associated to $ c $ has a quadratic relation if and only if $ t^2=1 $.
We calculate it for $ t=1 $ and $ a=p-2q $, and $ t=-1 $ and $ a=-p $.

\begin{pro} \label{prop:case1,7}
	Assume that  $ t=1 $ and $ a=p-2q $, or $ t=-1 $ and $ a=-p $. Then, the Nichols algebras are as in Table \ref{tab:h17}, where 
	\begin{align}\label{eqn:J_1,7_1}
	\lg  x_2x_1 - x_1x_2 - k \, x_1^2,   \, x_3x_1 - x_1x_3 - q \, x_1^2, \ \ \ \ \ \ \ \ \ \ \ \ \ \ \ \ \ \ \ \ \ \ \ \ \ \ \ \ \ \ \ \ \\
	\notag  x_3x_2 - x_2x_3 + k \, x_1x_3  -q \, x_1x_2  - (d+qk) \, x_1^2  \rg;  \\ \label{eqn:J_1,7_2}
	\lg  x_2x_1 + x_1x_2, \, x_3x_1 + x_1x_3, \, x_3x_2 + x_2x_3 - k \, x_1x_3 - q \, x_1x_2, \ \ \ \ \ \ \ \ \ \ \, \\
	\notag  x_1^2,\, x_2^2 - \ell \, x_1x_2 \rg; \\ \label{eqn:J_1,7_3}
	\lg  x_2x_1 + x_1x_2, \, x_3x_1 + x_1x_3, \, x_3x_2 +  x_2x_3 - k \, x_1x_3 - q \, x_1x_2, \ \ \ \ \ \ \ \ \ \ \, \\
	\notag x_1^2,\, x_2^2 - \ell \, x_1x_2, \, x_3^2 -p \, x_1x_3  \rg.
	\end{align}

	\begin{table}[ht]
		\caption{Nichols algebras of type $ \hit_{1, 7} $}\label{tab:h17}
		\begin{center}
			\begin{tabular}{| c | c | c | c | c | c | c |}\hline
				Case & $ t $ & $ a $ & $ b $ & $\J(V)$ &  Basis &$\GK$
				\\\hline
				\label{case1:17} \rm{(a)} & $1$ & $ p-2q $  &   & \eqref{eqn:J_1,7_1} & $ B_3 $ \eqref{eqn:B_3} & $ 3 $ \\\hline
				\label{case2:17} \rm{(b)} &$-1$ & $ -p $ & $ \neq -p^2 $   & \eqref{eqn:J_1,7_2} & $ B_1 $ \eqref{eqn:B_1} & $ 1 $ \\\hline
				\label{case3:17} \rm{(c)} &$-1$ & $ -p $ &  $  -p^2 $  & \eqref{eqn:J_1,7_3} & $ B_0 $ \eqref{eqn:B_0} & $ 0, \, \dim = 8 $ \\\hline
			\end{tabular}
		\end{center}
	\end{table}

\end{pro}

\pf  We get the quadratic relations through  \eqref{eq:J2}. 
Next we study the cases individually.

\smallbreak

\noindent Case \hyperref[case1:17]{\rm{(a)}}:
It follows similarly to Proposition \ref{prop:case1,5_t=1_and_p=2q}.
Some of the necessary identities to prove it are: 
\begin{align*}
&x_2x_1^{a_1} = x_1^{a_1}x_2 + a_1k \, x_1^{a_1+1},\qquad\quad x_3x_1^{a_1} = x_1^{a_1}x_3 + a_1q \, x_1^{a_1+1}, \\ 
x_3x_2^{a_2} = & \ x_2^{a_2}x_3 - a_2 k  x_1x_2^{a_2-1}x_3 +  a_2 q x_1 x_2^{a_2}  + (d+qk)\!  \sum_{\mathclap{i\in\I_{a_2}}}  \! \frac{a_2! \, k^{i-1}}{(a_2-i)!  i}  x_1^{i+1}x_2^{a_2-i}\!\!, \\
&\partial_2(x_1^{a_1}x_2^{a_2}x_3^{a_3}) = \sum_{\mathclap{i\in\I_{a_2}}} {a_2 \choose i} \ \prod_{\mathclap{j\in\I_{0, i-2}}} (jk+\ell) \, x_1^{a_1-1+i}x_2^{a_2-i}x_3^{a_3}, \\ 
& \ \ \ \   \partial_3(x_1^{a_1}x_3^{a_3}) = \sum_{\mathclap{i\in\I_{a_3}}} {a_3 \choose i} \ \prod_{\mathclap{j\in\I_{0, i-2}}} (p+jq) \, x_1^{a_1-1+i}x_3^{a_3-i}.
\end{align*}


\smallbreak

\noindent Case \hyperref[case2:17]{\rm{(b)}}:
Analogous to case \hyperref[case2:12]{\rm{(b)}} of Proposition \ref{prop:case1,2_remaining}.
We just present the derivations:
\begin{align*}
\partial_1(x_3^{a_3}) \! =& \left\lfloor a_3/2\right\rfloor \left( \left( \left\lfloor (a_3-2)/2\right\rfloor   p^2 -b \right)  x_1 + (-1)^{a_3} p \, x_3 \right) x_3^{a_3-2}, \\
\partial_1(x_1x_3^{a_3}) \! = & \, x_3^{a_3} + (-1)^{a_3 + 1} \left\lfloor {(a_3+1)}/{2}\right\rfloor p \, x_1x_3^{a_3-1}, \\ 
\partial_1(x_2x_3^{a_3}) \! = &  \left\lfloor{a_3}/{2}\right\rfloor \left( \left\lfloor {a_3}/{2}\right\rfloor  p^2 - b \right)  x_1x_2x_3^{a_3 - 2}  +  \chi_{_{o}}(a_3) \, (qk +d)  \, x_1x_3^{a_3 - 1} \\
& \,  +  (-1)^{a_3 + 1} \left\lfloor {(a_3 + 1)}/{2}\right\rfloor  p \,  (x_2 + k\, x_1)\, x_3^{a_3 - 1}, \\
\partial_1(x_1x_2x_3^{a_3}) \! = & \, (x_2 + k \, x_1) \, x_3^{a_3}   +  (-1)^{a_3} \left( \left\lfloor {(a_3 + 1)}/{2}\right\rfloor \! + \chi_{_{o}}(a_3)  \right) p \, x_1x_2x_3^{a_3 - 1} \\
\partial_2(x_1^{a_1}x_2^{a_2}x_3^{a_3}) \! =  & \, \delta_{a_2, 1} (-x_1)^{a_1}x_3^{a_3}, \\
\partial_3(x_1^{a_1}x_2^{a_2}x_3^{a_3}) \! = &  \, \chi_{_{o}}(a_3) (-x_1)^{a_1}(-k \, x_1 \! - x_2)^{a_2}x_3^{a_3-1} \! - \! \left\lfloor a_3/2\right\rfloor p \,  x_1^{a_1+1}x_2^{a_2}x_3^{a_3-2}.
\end{align*}

\smallbreak

\noindent Case \hyperref[case3:17]{\rm{(c)}}: Similar to the case $ t=-1 $ of Proposition \ref{prop:case1,6}.
\epf

\subsection{Case $ \hit_{1, 8}$}\label{subsec:hit18}
Let $ c $ be the braiding associated to the solution of QYBE  
$$ \hit_{1, 8} = \left(\begin{array}{ccc|ccc|ccc}
t & ta & tq & -ta & -ta^2 & -tb & -tq & tb & tap \\ 
\cdot & t & \cdot & \cdot & -ta & tq & \cdot & \cdot & tp \\ 
\cdot & \cdot & t & \cdot & \cdot & -ta & \cdot & \cdot & \cdot  \\ \hline
\cdot & \cdot & \cdot & t & ta & \cdot & \cdot & -tq & -tp \\ 
\cdot & \cdot & \cdot & \cdot & t & \cdot & \cdot & \cdot & \cdot \\ 
\cdot & \cdot & \cdot & \cdot & \cdot & t & \cdot & \cdot & \cdot \\ \hline
\cdot & \cdot & \cdot & \cdot & \cdot & \cdot & t & ta & \cdot  \\ 
\cdot & \cdot & \cdot & \cdot & \cdot & \cdot & \cdot & t & \cdot  \\ 
\cdot & \cdot & \cdot & \cdot & \cdot & \cdot & \cdot & \cdot & t 
\end{array} \right). $$
For the presentation of $ c $, utilize \eqref{eqn:QYBE_to_braiding}.


\begin{pro} \label{prop:case1,8}
	If $ t^2 \neq 1 $, then there are no quadratic relations. Otherwise, the Nichols algebras are as in Table \ref{tab:h18}, where 
	\begin{align}\label{eqn:Nichols1,8_t=-1}
	\lg  x_2x_1 + x_1x_2, \, x_3x_1 + x_1x_3, & \ x_3x_2 +  x_2x_3 - a \, x_1x_3 + q \, x_1x_2, \\
	\notag  & \qquad \qquad \  x_1^2,\, x_2^2 - a\, x_1x_2, \, x_3^2 +p\, x_1x_2  \rg; \\
	\label{eqn:Nichols1,8_t=1}
	\lg  x_2x_1 - x_1x_2 - a \, x_1^2,  \, x_3x_1 - \, &  x_1x_3 - q \, x_1^2,  \\
	\notag  & \, x_3x_2 - x_2x_3 + a \, x_1x_3 - q \, x_1x_2 + b \, x_1^2 \rg.
	\end{align}
	
	\begin{table}[ht]
		\caption{Nichols algebras of type $ \hit_{1, 8} $}\label{tab:h18}
		\begin{center}
			\begin{tabular}{| c | c | c | c | c | c |}\hline
				$ t $  & $\J(V)$ &  Basis &$\GK$
				\\\hline
				$-1$ &  \eqref{eqn:Nichols1,8_t=-1} & $ B_0 $ \eqref{eqn:B_0} & $ 0, \, \dim = 8 $ \\\hline
				$1$ &  \eqref{eqn:Nichols1,8_t=1} & $ B_3  $ \eqref{eqn:B_3} & $ 3 $ \\\hline
			\end{tabular}
		\end{center}
	\end{table}	
\end{pro}

\pf Analogous to Proposition \ref{prop:case1,6}.
Some of the necessary equations to prove case $ t=1 $ are: 
\begin{align*}
x_2x_1^{a_1} = x_1^{a_1}x_2 + a_1 a \, x_1^{a_1+1},   & \qquad\qquad x_3x_1^{a_1} = x_1^{a_1}x_3 + a_1q \, x_1^{a_1+1}, \\ 
 x_3x_2^{a_2} = x_2^{a_2}x_3 - a_2 a \, x_1x_2^{a_2-1}& x_3  +  a_2 q \, x_1 x_2^{a_2}    - b   \sum_{\mathclap{i\in\I_{a_2}}} \frac{a_2! \, a^{i-1}}{(a_2-i)! \, i} \, x_1^{i+1}x_2^{a_2-i}, \\
\partial_3(x_1^{a_1}x_2^{a_2}x_3^{a_3}) = & \, a_3 \sum_{\mathclap{i\in\I_{0, a_2}}} \frac{a_2 ! \,a^i }{(a_2 - i)!} \, x_1^{a_1+i}x_2^{a_2-i}x_3^{a_3-1}, \\ 
\partial_2(x_1^{a_1}x_2^{a_2}) = \, & \sum_{\mathclap{i\in\I_{a_2}}} \frac{a_2 ! \, a^{i-1} }{(a_2 - i)! \, i} \, x_1^{a_1-1+i}x_2^{a_2-i}.
\end{align*}
\epf

\subsection{Case $ \hit_{1, 9}$}\label{subsec:hit19}
Let $ c $ be the braiding associated to the solution of QYBE  
$$ \hit_{1, 9} = \left(\begin{array}{ccc|ccc|ccc}
t & t & \cdot & -t & -t & ta & tb & t(b-a) & tp \\ 
\cdot & t & t & \cdot & -t & -t & \cdot & tb & t(b-q) \\ 
\cdot & \cdot & t & \cdot & \cdot & -t & \cdot & \cdot & tb  \\ \hline
\cdot & \cdot & \cdot & t & t & \cdot & -t & -t & tq \\ 
\cdot & \cdot & \cdot & \cdot & t & t & \cdot & -t & -t \\ 
\cdot & \cdot & \cdot & \cdot & \cdot & t & \cdot & \cdot & -t \\ \hline
\cdot & \cdot & \cdot & \cdot & \cdot & \cdot & t & t & \cdot  \\ 
\cdot & \cdot & \cdot & \cdot & \cdot & \cdot & \cdot & t & t  \\ 
\cdot & \cdot & \cdot & \cdot & \cdot & \cdot & \cdot & \cdot & t 
\end{array} \right). $$
Apply \eqref{eqn:QYBE_to_braiding} for the presentation of $ c $.


As previously, the Nichols algebra associated to $ c $ has a quadratic relation iff $ t^2=1 $.
We compute it just when $ b=1 $ also.

\begin{pro} \label{prop:case1,9}
	\noindent Assume that $ t^2 \neq 1 $ and $ b=1 $. Then the Nichols algebras are as in Table \ref{tab:h19}, where 
	\begin{align}\label{eqn:J_1,9_1}
	\lg  x_2x_1 - x_1x_2 - x_1^2, \,   x_3 & x_1 - x_1x_3 - x_1x_2, \, \\
	\notag & \hspace{0.8cm} x_3x_2 - x_2x_3 - x_2^2  + x_1x_3  + x_1x_2  -a \, x_1^2  \rg;  \\ \label{eqn:J_1,9_4}
	\lg  x_2x_1 + x_1x_2, \, x_3x_1 + & \, x_1x_3 + x_1x_2, \, x_3x_2 +  x_2x_3 -   x_1x_3, \, \\
	\notag & \hspace{4.6cm} x_1^2,\, x_2^2 - x_1x_2 \rg; \\ \label{eqn:J_1,9_5}
	\lg  x_2x_1 + x_1x_2, \, x_3x_1 +  & \, x_1x_3 + x_1x_2,  \, x_3x_2 + x_2x_3 - x_1x_3,  \\
	\notag  & \, x_1^2,\, x_2^2 - x_1x_2, \, x_3^2 - x_2x_3 + x_1x_3 - q\, x_1x_2  \rg.
	\end{align}
	
	\begin{table}[ht]
		\caption{Nichols algebras of type $ \hit_{1, 9} $}\label{tab:h19}
		\begin{center}
			\begin{tabular}{| c | c | c | c | c | c | c |}\hline
			Case &	$ t $ & $ b $ & $ p $ & $\J(V)$ &  Basis &$\GK$
				\\\hline
				\label{case1:19} \rm{(a)} &$1$ & $ 1 $  &   & \eqref{eqn:J_1,9_1} & $ B_3 $ \eqref{eqn:B_3} & $ 3 $ \\\hline
				\label{case2:19} \rm{(b)} &$-1$ & $ 1 $ & $ \neq -a-q $   & \eqref{eqn:J_1,9_4} & $ B_1 $ \eqref{eqn:B_1} & $ 1 $ \\\hline
				\label{case3:19} \rm{(c)} &$-1$ & $ 1 $ &  $  -a-q $  & \eqref{eqn:J_1,9_5} & $ B_0 $ \eqref{eqn:B_0} & $ 0, \, \dim = 8 $ \\\hline
			\end{tabular}
		\end{center}
	\end{table}	 
	
\end{pro}

\pf  

The quadratic relations above hold by \eqref{eq:J2}. 
Next we proceed the analysis of each case enumerated above separately.

\smallbreak

\noindent Case \hyperref[case1:19]{\rm{(a)}}:
We proceed as in case $ t=1 $ of Proposition \ref{prop:case1,1}. First, we prove
\begin{align*}
x_3x_1^{a_1} =  x_1^{a_1}(x_3 &+ a_1  x_2) + \binom{a_1}{2} \, x_1^{a_1+1}, \qquad   x_2^{ a_2}x_1^{a_1} = \sum_{\mathclap{i\in\I_{0, a_2}}} T_{a_1,a_2}^{i}  x_1^{a_1+i}x_2^{a_2-i}, \\ 
x_3x_2^{a_2} &= x_2^{a_2}x_3 + a_2 \, x_2^{a_2+1} - a_2 \, x_1x_2^{a_2-1}x_3   \\
 & \ \ \ - {a_2+1\choose 2} \, x_1 x_2^{a_2} + a \,  \sum_{\mathclap{i\in\I_{a_2}}} \frac{a_2!}{(a_2-i)! \, i} \, x_1^{i+1}x_2^{a_2-i},
\end{align*}
for all $ a_1, a_2 \geq 0 $, where  $T_{a_1,a_2}^{i} = \frac{(a_1+i-1) !}{(a_1 -1)!}{{a_2}\choose {i}}$ if $ a_1\geq 1$, and $ T_{0,a_2}^{i} = \delta_{i, 0}$. 
We also have
\begin{align*}
\partial_3(x_1^{a_1}x_2^{a_2}x_3^{a_3}) &= \sum_{\mathclap{i\in\I_{0, a_2}}} \frac{a_2!}{(a_2-i)!} \, x_1^{a_1+i}x_2^{a_2-i}\partial_3 (x_3^{a_3}),
\end{align*}
what implies that
\begin{align*}
\notag 0 &= \! \! \! \! \! \! \! \! \! \! \! \sum_{ \substack{a_1+a_2+a_3=n \\ a_i\geq 0}}  \! \! \! \! \! \! \! \! \!  \lambda_{a_1, a_2, a_3}  \! \!  \!  \! \! \! \! \! \! \! \! \!   \sum_{\substack{c_1+c_2+c_3=a_3-1 \\ c_i\geq 0}} \! \!  \!  \! \! \! \! \! \! \! \! \! \!  \mu_{a_3, c_1, c_2, c_3} \!\!\! \sum_{i\in\I_{0, a_2}}\sum_{j\in\I_{0, a_2-i}}  \!\!\!\!\!\! \dfrac{a_2! \, T^j_{c_1, a_2-i}}{(a_2-i)!}  x_1^{a_1+c_1+i+j}\! x_2^{a_2+c_2-i-j}  x_3^{c_3}\! .
\end{align*}

Finally we vanish all $ \lambda_{a_i} $ as in Proposition \ref{prop:case1,5_t=1_and_p=2q}. 
To do so, we also use that 
\begin{align*}
\partial_2(x_1^{a_1}x_2^{a_2}) &= \sum_{\mathclap{i\in\I_{a_2}}} \frac{a_2!}{(a_2-i)! \, i}\, x_1^{a_1-1+i}x_2^{a_2-i}.
\end{align*}

\smallbreak

\noindent Case \hyperref[case2:19]{\rm{(b)}}: 
Similar to case \hyperref[case2:12]{\rm{(b)}} of Proposition \ref{prop:case1,2_remaining}.
We only present the derivations.
\begin{align*}
&\partial_1(x_3^{a_3})  =  (-1)^{a_3+1} \left\lfloor \dfrac{a_3}{2}\right\rfloor \left\lfloor \dfrac{a_3+1}{2}\right\rfloor
x_3^{a_3-1} -( p + \chi_{_{o}}(a_3) \, a) x_1x_3^{a_3-2} \\
&\hspace{0.62cm} +\left\lfloor \dfrac{a_3}{2}\right\rfloor \! q \, x_2x_3^{a_3-2} \! + \! \left\lfloor \dfrac{a_3}{2}\right\rfloor \! \left(  \left\lfloor \dfrac{a_3}{2}\right\rfloor^3 \!\! +  \chi_{_{o}}(a_3)  \dfrac{a_3(a_3-1)}{2} \right) \! x_1(x_2+x_3)x_3^{a_3-3}  \\
&\hspace{0.62cm} - \left\lfloor \dfrac{a_3}{2}\right\rfloor \!\! \left( \dfrac{2}{3}\left\lfloor \dfrac{a_3}{2}\right\rfloor^2 \!\! + \dfrac{1}{2} \left\lfloor \dfrac{a_3}{2}\right\rfloor \! - \dfrac{1}{6} \right) \! (x_1x_3+x_2x_3+x_1x_2)x_3^{a_3-3} \\
&\hspace{0.62cm}  - \left( \! \chi_{_{e}}(a_3) \dfrac{(a_3-2)(q+a)}{2} + \chi_{_{o}}(a_3)(p+a)\! \right) \! x_1x_2x_3^{a_3-3}, \\
&\partial_1(x_2x_3^{a_3})   =  (-1)^{a_3} \left\lfloor \dfrac{a_3+1}{2}\right\rfloor   \left\lfloor \dfrac{a_3+2}{2}\right\rfloor (x_1+x_2)x_3^{a_3-1}  \\
&\hspace{0.62cm} + \!\left\lfloor \dfrac{a_3}{2}\right\rfloor \!  \left(  \left\lfloor \dfrac{a_3}{2}\right\rfloor^3 \! + 2  \left\lfloor \dfrac{a_3}{2}\right\rfloor^2 \! + \left\lfloor \dfrac{a_3}{2}\right\rfloor + \left\lfloor \dfrac{a_3+3}{2}\right\rfloor  - p - q - a \! \right) \! x_1x_2x_3^{a_3-2}\\
&\hspace{0.62cm} + \! \chi_{_{o}}(a_3) \left\lfloor \dfrac{a_3}{2}\right\rfloor  \left( 2 \left\lfloor \dfrac{a_3}{2}\right\rfloor^2 \! + 3  \left\lfloor \dfrac{a_3}{2}\right\rfloor \! - a \right) x_1x_2x_3^{a_3-2} + \chi_{_{o}}(a_3) \, a \, x_1x_3^{a_3-1} \! , \\ 
&\partial_1(x_1x_2^{a_2}x_3^{a_3})  =  x_1^{a_2}x_3^{a_3} +  \left\lfloor \dfrac{a_3}{2}\right\rfloor  \left(  \dfrac{2}{3}\left\lfloor \dfrac{a_3}{2}\right\rfloor^2 + \dfrac{3}{2}\left\lfloor \dfrac{a_3}{2}\right\rfloor  + \dfrac{5}{6}  - q \right) x_1x_2^{a_2}x_3^{a_3-2} \\
&\hspace{0.62cm} + \! \left( \! (-1)^{a_3}  \! \left\lfloor \dfrac{a_3+1}{2}\right\rfloor \!\!\left\lfloor \dfrac{a_3+2+2a_2}{2}\right\rfloor \! -a_2 \chi_{_{o}}(a_3) \! \right) \! x_1(-x_2)^{a_2}x_3^{a_3-1} \! + x_2^{a_2}x_3^{a_3} \! , \\
&\partial_2(x_1^{a_1}x_2^{a_2}x_3^{a_3})  =   \delta_{a_2, 1} (-x_1)^{a_1}x_3^{a_3}  -  \left\lfloor a_3/2\right\rfloor q\,    x_1^{a_1+1}x_2^{a_2}(\chi_{_{o}}(a_3)\, x_2+x_3)x_3^{a_3-3} \\
&\hspace{0.3cm} -\! \left\lfloor \dfrac{a_3}{2}\right\rfloor \!\! \left(  \dfrac{4}{3}\left\lfloor \dfrac{a_3}{2}\right\rfloor^2 \!\!\! + \! \left( \! 2\chi_{_{o}}(a_3) + \dfrac{4\, a_2-3}{2}\right) \!\! \left\lfloor \dfrac{a_3}{2}\right\rfloor \! + \dfrac{1}{6}   \right) \! x_1^{a_1+1}x_2^{a_2}(x_2+x_3)x_3^{a_3-3} \\
& \hspace{0.3cm} + \! (-1)^{a_3}(\left\lfloor (a_3+a_1+a_2)/2\right\rfloor+a_1a_2(-1)^{a_3-1})   (-x_1)^{a_1}(-x_1-x_2)^{a_2}x_3^{a_3-1} \\
& \hspace{0.3cm} + \! \left\lfloor a_3/2\right\rfloor \left\lfloor (a_3+2a_1+3a_2)/2\right\rfloor   (-x_1)^{a_1}(-x_2)^{a_2}x_2x_3^{a_3-2}, \\
&\partial_3(x_1^{a_1}x_2^{a_2}x_3^{a_3})  =  \left\lfloor a_3/2\right\rfloor \left\lfloor (a_3-1)/2\right\rfloor   x_1^{a_1+1}x_2^{a_2}(x_2+x_3)x_3^{a_3-3} \\
& +  \chi_{_{o}}(a_3) (-x_1)^{a_1}(- x_1 \! - x_2)^{a_2}x_3^{a_3-1} \! + (a_2+1) \! \left\lfloor a_3/2\right\rfloor \!   (-x_1)^{a_1}(-x_2)^{a_2+1}x_3^{a_3-2}\! .
\end{align*}

\smallbreak

\noindent Case \hyperref[case3:19]{\rm{(c)}}: Analogous to case $ t=-1 $ of Proposition \ref{prop:case1,6}.
\epf

\subsection{Case $ \hit_{1, 10}$}\label{subsec:hit110}
Let $ c $ be the braiding associated to the solution of QYBE $ \hit_{1, 10} = $ 
{\footnotesize $$  t\left(\begin{array}{ccc|ccc|ccc}
	1 & 2 & \cdot & -2 & -4a & 8a & 4 & \cdot & -4a(a+7) \\ 
	\cdot & 1 & 2 & \cdot & -2a & 4-6a-2a^2 & \cdot & 2a(a+1) & 6a(a+2)(a-1) \\ 
	\cdot & \cdot & 1 & \cdot & \cdot & 2-4a & \cdot & \cdot & -4a(1-2a)  \\ \hline
	\cdot & \cdot & \cdot & 1 & 2a & -2a(1-a) & -2 & -2a(a+1) & -2a(1+3a)(a-2) \\ 
	\cdot & \cdot & \cdot & \cdot & 1 & 2a & \cdot & -2a & 4a(1-2a) \\ 
	\cdot & \cdot & \cdot & \cdot & \cdot & 1 & \cdot & \cdot & 2-4a \\ \hline
	\cdot & \cdot & \cdot & \cdot & \cdot & \cdot & 1 & 4a-2 & 4(1-2a)(1-a)  \\ 
	\cdot & \cdot & \cdot & \cdot & \cdot & \cdot & \cdot & 1 & 4a-2  \\ 
	\cdot & \cdot & \cdot & \cdot & \cdot & \cdot & \cdot & \cdot & 1 
	\end{array} \right). $$}
Utilize \eqref{eqn:QYBE_to_braiding} for the explicit presentation of $ c $.


\begin{pro} \label{prop:case1,10}
	If $ t^2 \neq 1 $, then there are no quadratic relations. Otherwise, the Nichols algebras are as in Table \ref{tab:h110}, where 
	\begin{align}\label{eqn:J_1,10_1}
	\lg &  x_2 x_1 - x_1 x_2 - 2 \, x_1^2,  \, x_3 x_1 - x_1 x_3 - 2 \, x_1 x_2, x_3 x_2 -  x_2 x_3 \\
	\notag		& \, \ \ \ \ \ \ \ \ \ \  - 2 a \, x_2^2 + (4a-2) x_1 x_3  + (8a-4) x_1 x_2 - 4 a (a+1) \, x_1^2  \rg;\\ \label{eqn:J_1,10_2}
	\lg &   x_2x_1 + x_1x_2,  \, x_3x_1 + x_1x_3 + 2 \, x_1x_2, \, x_3x_2 + x_2x_3  \\
	\notag \, & \ \ \ \ \ \ \ \ \ \ \ \ \ \ \ \ + 2(1-2a)\, x_1x_3 + 4(1-a)\, x_1x_2, \, x_1^2, \, x_2^2 - 2a \, x_1x_2 \rg; \\ \label{eqn:J_1,10_3}
	\lg &  x_2x_1 +  x_1x_2, \, x_3x_1 + x_1x_3 + 2 \, x_1x_2, \, x_3x_2 + x_2x_3 + 2(1-2a)\, x_1x_3 \\
	\notag  & \ \ \ \ \ \ \ \ \ \ \ \    + 4(1-a)\, x_1x_2, \, x_1^2, \, x_2^2 - 2a \, x_1x_2, \, x_3^2 -2(2a-1)\, x_2x_3 \\
	\notag  & \ \ \ \ \ \ \ \ \ \ \ \ \ \ \ \ \ \ \ \ \ \ \ \   - a(a^2-2a-3)\, x_1x_3 -2a(a^2 - a+4)\, x_1x_2  \rg.
	\end{align}

	\begin{table}[ht]
		\caption{Nichols algebras of type $ \hit_{1, 10} $}\label{tab:h110}
		\begin{center}
			\begin{tabular}{| c | c | c | c | c | c |}\hline
				Case & $ t $ & $ a $  & $\J(V)$ &  Basis &$\GK$
				\\\hline
				\label{case1:110} \rm{(a)} &$1$ &      & \eqref{eqn:J_1,10_1} & $ B_3 $ \eqref{eqn:B_3} & $ 3 $ \\\hline
				\label{case2:110} \rm{(b)} &$-1$ & $\notin\{-7, \, 0, \, 1\}$    & \eqref{eqn:J_1,10_2} & $ B_1 $ \eqref{eqn:B_1} & $ 1 $ \\\hline
				\label{case3:110} \rm{(c)} &$-1$ & $\in\{-7, \, 0, \, 1\}$     & \eqref{eqn:J_1,10_3} & $ B_0 $ \eqref{eqn:B_0}  & $ 0, \, \dim = 8 $ \\\hline
			\end{tabular}
		\end{center}
	\end{table}

\end{pro}

\pf The relations above are obtained via  \eqref{eq:J2}. 
Next we study the cases individually.

\smallbreak

\noindent Case \hyperref[case1:110]{\rm{(a)}}: This case follows similarly to the case $(a)$ of Proposition \ref{prop:case1,9}.
We have that
\begin{align*}
x_3x_1^{a_1} =& \ x_1^{a_1}\left( \! x_3 + 2 a_1  x_2 + 4  \binom{a_1}{2}  x_1 \! \right),  
\qquad x_2^{ a_2}x_1^{a_1} = \sum_{\mathclap{i\in\I_{0, a_2}}} 2^i  T_{a_1,a_2}^{i}  x_1^{a_1+i}x_2^{a_2-i},  \\
x_3x_2^{a_2} =& \ x_2^{a_2}x_3 + 2\, a \, a_2 \, x_2^{a_2+1} +   \sum_{\mathclap{i\in\I_{a_2}}} A^{a_2}_i\, x_1^{i+1}x_2^{a_2-i} + \sum_{\mathclap{i\in\I_{a_2-1}}} B^{a_2}_i\, x_1^{i+1}x_2^{a_2-1-i}x_3\\
& - a_2 (4 a-2) \, x_1x_2^{a_2-1}x_3 - a_2(2a-1)(2a (a_2-1)+4) \, x_1 x_2^{a_2},  
\end{align*}
for some $ A_i^{a_2}, B_i^{a_2}\in \ku $ and  $ T_{a_1,a_2}^{i} $ as in Proposition \ref{prop:case1,9}.
It also holds that
\begin{align*}
\partial_3(x_1^{a_1}x_2^{a_2}x_3^{a_3}) &= \sum_{\mathclap{i\in \I_{0,a_2}}}  2^{\left\lfloor \frac{3i}{2} \right\rfloor} { a_2 \choose i} \, C_i \, D_{\left\lfloor \frac{i-2}{2} \right\rfloor} \, x_1^{a_1+i}x_2^{a_2-i} \partial_3 (x_3^{a_3}), \\
\partial_2(x_1^{a_1}x_2^{a_2}) &= \sum_{\mathclap{i\in \I_{a_2}}} 2^{i-1} {a_2 \choose i} \, D_{i-2}  \, x_1^{a_1-1+i}x_2^{a_2-i},
\end{align*}
where $ C_i = \prod_{j\in\I_{0,\left\lfloor \frac{i-1}{2} \right\rfloor}} (2a+2j-1) $ and $ D_i = \prod_{j\in\I_{0,i}} (a+j) $.

\smallbreak

\noindent Case \hyperref[case1:110]{\rm{(b)}}: We proceed similarly to the case \hyperref[case2:12]{\rm{(b)}} of Proposition \ref{prop:case1,2_remaining}.
	We just present the derivations of the term $x_3^{a_3}$.
\begin{align*}
\partial_1(x_3^{a_3}) &= 
(-1)^{a_3+1}  4(2a-1) \left\lfloor {a_3}/{2} \right\rfloor \left( (2a-1) \left\lfloor {(a_3 - 1)}/{2} \right\rfloor + a \right) \  x_3^{a_3-1}\\
&\hspace{-7pt} - {2}/{3}\left\lfloor {a_3}/{2} \right\rfloor \left(E^{1,\, a_3}_{\left\lfloor {a_3}/{2} \right\rfloor} \, x_2x_3 - E^{2,\, a_3}_{\left\lfloor {a_3}/{2} \right\rfloor}  \, x_1x_3 + E^{3, \, a_3}_{\left\lfloor {a_3}/{2} \right\rfloor} \, x_1x_2 \right) x_3^{a_3-3},\\ 
\partial_2(x_3^{a_3})& = 2(2a-1)\! \left\lfloor {a_3}/{2} \right\rfloor \! \left( (-1)^{a_3} x_3 + 
\left((a_3 -1 )(2a-1) + \chi_{_{e}} (a_3) \right) x_2 \right) x_3^{a_3-2} \\
&\hspace{-7pt} - {2}/{3} \left\lfloor {a_3}/{2} \right\rfloor \left(F^{1, \,a_3}_{\left\lfloor {a_3}/{2} \right\rfloor} \,  x_1x_3^{a_3-2} +  F^{2,\, a_3}_{\left\lfloor {a_3}/{2} \right\rfloor}  \, x_1x_2x_3^{a_3-3}\right),\\  
\partial_3(x_3^{a_3})& =  4 a (2a-1)^2 \left\lfloor {(a_3 - 1)^2}/{2} \right\rfloor  x_1x_2x_3^{a_3-3} - 2(2a-1) \left\lfloor {a_3}/{2} \right\rfloor    x_2x_3^{a_3-2}\\
&\hspace{-7pt} + 4(2a-1) \left\lfloor {a_3}/{2} \right\rfloor \left( (2a-1) \left\lfloor {(a_3 - 3)}/{2} \right\rfloor +a  \right) x_1x_3^{a_3-2} +\chi_{_{o}}(a_3) \, x_3^{a_3 - 1}, 
\end{align*}
where
\begin{align*}
E^{1,n}_{r} &=  a^3(64r^2 + 24(1 - 4\chi_{_{e}}(n))r + (5 + 12\chi_{_{e}}(n)))\\
&- 3a^2 (32 r^2 + 8(1 - 8\chi_{_{e}}(n))r + 7(1 + 2\chi_{_{e}}(n))) \\
&+ 6a(8r^2 + (1 - 20\chi_{_{e}}(n))r + (1 + 7\chi_{_{e}}(n))) \\
&- 4(2r^2 - 6\chi_{_{e}}(n)r + (1 + 3\chi_{_{e}}(n))), \\
E^{2,n}_{r} &=  96a^4(4r^3  - 8\chi_{_{e}}(n)r^2 - (1 - 6\chi_{_{e}}(n))r - \chi_{_{e}}(n))\\
&- 16a^3(48r^3 - 16(1 + 6\chi_{_{e}}(n))r^2 - 3(5 - 34\chi_{_{e}}(n))r + (1 - 24\chi_{_{e}}(n))) \\
&+ 6a^2(96r^3 - 64(1 + 3\chi_{_{e}}(n))r^2 - 4(7 - 66\chi_{_{e}}(n))r + (1 - 76\chi_{_{e}}(n))) \\
&- 6a(32 r^3 -32(1 + 2\chi_{_{e}}(n))r^2 - 6(1 - 18\chi_{_{e}}(n))r - 5(1 + 8\chi_{_{e}}(n)))\\
&+ 8(3r^3 - 2(2 + 3\chi_{_{e}}(n))r^2 + 12\chi_{_{e}}(n)r + (1 - 6\chi_{_{e}}(n))), \\
E^{3,n}_{r} &=  2a^5( 64r^2 - 12(7 + 2\chi_{_{e}}(n)) r + (5 +39\chi_{_{e}}(n))) \\
&- 4a^4(192 r^3 \! + 16(7 - 24\chi_{_{e}}(n))r^2 - 27(9 - 10\chi_{_{e}}(n))r \! + (29 + 24\chi_{_{e}}(n)))\\
&+ 6a^3 (256 r^3 \! + \! 16(1 - 32\chi_{_{e}}(n))r^2 \! - \! 2(125 - \! 258\chi_{_{e}}(n))r \! + \!(29 - 55\chi_{_{e}}(n))) \\
&- 8a^2(144 r^3 \! - \! 2(23 + 144\chi_{_{e}}(n))r^2 \! - \! 3(38 - 115\chi_{_{e}}(n))r \! - (8 + 33\chi_{_{e}}(n))) \\
&+ 4a (96 r^3 -64(1 + 3\chi_{_{e}}(n))r^2 \! - 6(12 - 47\chi_{_{e}}(n))r - 5(7 + 3\chi_{_{e}}(n))) \\
&- 48( r^3 - (1 + 2\chi_{_{e}}(n)) r^2 - (1 - 4\chi_{_{e}}(n)) r + (1 - 2\chi_{_{e}}(n))), \\
F^{1,n}_{r} &=  a^3(128r^2 - 24(1+ 8\chi_{_{e}}(n)) r - (5 - 84\chi_{_{e}}(n)))\\
&- 3a^2 (64 r^2 - 24(1 + 4\chi_{_{e}}(n))r - (7 - 58\chi_{_{e}}(n)))\\
&+ 6a(16 r^2 - 3(3 + 8\chi_{_{e}}(n))r - (1 - 19\chi_{_{e}}(n))) \\
&-4(4r^2 - 3(1 + 2\chi_{_{e}}(n)) r - (1 - 6\chi_{_{e}}(n))), \\
F^{2,n}_{r} &=  2a^4 (64 r^2 - 24(1 + 4\chi_{_{e}}(n)) r + (5 + 51\chi_{_{e}}(n)))\\
&- 2a^3(32 r^2 - 24(3 + 2\chi_{_{e}}(n)) r + (31 + 57\chi_{_{e}}(n))) \\
&- 12a^2(8r^2 + 3(3 - 4\chi_{_{e}}(n)) r - (8 - 3\chi_{_{e}}(n))) \\
&+ 8a(10 r^2 + 3(1 - \chi_{_{e}}(n)5) r - 2(2 - 3\chi_{_{e}}(n))) \\
&- 8(2r^2 - 3\chi_{_{e}}(n)r - (2 - 3\chi_{_{e}}(n))).
\end{align*}

\smallbreak

\noindent Case \hyperref[case1:110]{\rm{(c)}}: Similar to case $ t=-1 $ of Proposition \ref{prop:case1,6}.
\epf

{\color{magenta}

%
%
%
%
%
%
%

}

\end{document}